\documentclass[12pt]{elsarticle}
\usepackage{amsmath,amsfonts,amssymb,graphicx,
verbatim}

 \newcommand{\bff}[1]{{\mbox{\boldmath$#1$}}}

\def\mP{{\mathcal P}}
  \def\mM{{\mathcal M}}

\def\rui{\Rui}
\def\srui{\ovl{\rui}}
\def\Rui{\Psi}  
 \def\Rui{\Psi}

\def\BEN{\begin{enumerate}}  \def\BI{\begin{itemize}}
\def\EEN{\end{enumerate}}   \def\EI{\end{itemize}}
\def\beq{\begin{eqnarray}} \def\eeq{\end{eqnarray}}
\def\bea{\begin{eqnarray*}}
\def\eea{\end{eqnarray*}}
\def\le{\left} \def\ri{\right} 
\def\te#1{\mathrm{e}^{#1}}
\def\T{\tilde}  \def\H{\hat}
\def\ovl{\overline}

\def\I{\infty}
\def\a{\alpha} \def\b{\beta}
\def\g{\gamma}  \def\d{\delta}   \def\th{\theta}
 \def\k{\kappa} \def\l{\lambda} \def\m{s} \def\n{\eta}
  \def\r{\rho} \def\s{\sigma}
 \def\p{\psi}  \def\f{\psi^*}  
\def\c{c} \def\w{\omega} \def\q{\qquad} \def\D{\Delta}
 \def\G{\Gamma}

\def\la{\label} \def\fr{\frac} \def\im{\item}

\newtheorem{thm}{Theorem}
\def\beT{\begin{thm}
  }
  \def\eeT{\end{thm}}

  \newtheorem{Qu}{Question}
\def\beQ{\begin{Qu}}
  \def\eeQ{\end{Qu}}

\newtheorem{Lem}{Lemma}
\newtheorem{ex}{Example}
\newtheorem{Pro}{Proposition}

\newtheorem{Def}{Definition}

\newtheorem{rem}{Remark}

\newtheorem{Exe}{Exercice}
\def\beXe{\begin{Exe}} \def\eeXe{\end{Exe}}

\def\eeD{\end{Def}} \def\beD{\begin{Def}}
\def\beXa{\begin{ex}} \def\eeXa{\end{ex}}
\def\beR{\begin{rem}} \def\eeR{\end{rem}}
\def\beL{\begin{Lem}} \def\eeL{\end{Lem}}
\def\beP{\begin{Pro}} \def\eeP{\end{Pro}}
\def\beC{}

\long\def\symbolfootnote[#1]#2{
\begingroup
\def\thefootnote{\fnsymbol{footnote}}\footnote[#1]{#2}
\endgroup}
\def\fn{\symbolfootnote}
 \def\sec{\section}
\def\WT{\widetilde}
\def\WH{\widehat}
\def\no{\nonumber} 
\def\bep{\begin{pmatrix}}
  \def\eep{\end{pmatrix}}
  \def\Eq{\Leftrightarrow}
  \def\CL{Cram\'er-Lundberg }
  \def\PK{Pollaczek Khinchine }

    \def\E{{\mathbb E}}

\begin{document}
\title{
 On
 matrix exponential
 approximations\\
of the infimum of a \\spectrally negative Levy process}

 \author[uppa]{F. Avram}
 \author[dit]{A. Horv\'ath}
 \author[uy1]{M. R. Pistorius}
 \cortext[cor1]{Corresponding author, Tel  : +33 (0)5 59 40 75 63, Fax : +33 (0)5 59 40 75 55}
\address[uppa]{ Labo. de Math. Appliqu\'ees,
Universit\'e de Pau et des Pays de l'Adour, 
Pau,  France}
\address[dit]{Dipart. di Informatica. Universit\`a di
Torino Corso Svizzera, 185 10149, Torino, Italy}
\address[uy1]{Department of Mathematics Imperial
College,  London, United Kingdom}

\begin{abstract}
We recall four open problems concerning constructing high-order
 {\bf matrix-exponential approximations} for the 
 infimum of a spectrally negative Levy process (with applications to {\bf first-passage/ruin
probabilities},
  the {\bf waiting time distribution} in the M/G/1
queue,  pricing of barrier options, etc).

On the way, we provide   a new approximation, for the perturbed \CL model, and recall  a remarkable family of   (not minimal order) approximations of Johnson and
Taaffe \cite{johnson1989matching}, which   fit an {\bf arbitrarily high  number of moments},
   greatly generalizing the  currently used %
approximations of Renyi, De Vylder and Whitt-Ramsay.  Obtaining such approximations which fit   the Laplace transform at infinity as well would be quite useful.
\end{abstract}

\maketitle

 {\bf keywords:} {Levy process; first passage problem;
 Pollaczek-Khinchine formula;  method of moments;
 matrix-exponential function;
 admissible Pad\'e approximation; 
 Johnson-Taaffe approximations;
  two-point Pad\'e approximations}



%

\section{Introduction}
\vskip.1truecm{\bf Motivation:}
The  problem of
approximating  distributions   based on  empirical data like moments
is one of the bread and butter problems of applied probability.

In
  risk theory and related {\bf first passage
  applications}  (queueing,
mathematical finance, ...) there is special interest in
approximating
densities of nonnegative random variables by {\bf affine combinations of exponentials}, also
called {\bf GHE densities (generalized hyper-exponential)}, or, more generally, by  {\bf matrix exponential distributions} (which allow for "collision of exponents"). One reason for that is that this class 
captures the asymptotic behavior in the important "light tails case".

There exists  already a quite extensive literature, based on inverting
the explicit Pollaczek-Khinchine formula for the Laplace transform, which assumes
  however  complete knowledge of the model, for example knowledge of all its moments.

  Since  input data are never certain, it is interesting to  develop {\bf approximations based on finitely many moments} (the  coefficients of the power series expansion of the Laplace transform
around $0$).

Some well-known 
such approximations in risk theory and queueing   are  the    {\bf Renyi, De Vylder, Gamma},  and {\bf Whitt -Ramsay} approximations,
obtained by {\bf fitting one, two or three moments}, and the minimal (but arbitrarily high) order {\bf three moments fitting formulas of Bobbio, Horvath and Telek}
\cite{bobbio2005matching}.
Also useful    is the {\bf Cram\'er-Lundberg}
approximation \beq \la{CLa}\rui(x) \sim C e^{- \g x}\eeq
 where $\rui(x)$ denotes the ruin probability \eqref{eq:ru} and $-\g$ is the so called adjustment coefficient (i.e. the
largest negative root of the Cram\'er-Lundberg equation \eqref{eq:CL}), which captures  the asymptotic
behavior in the case of light tail claims.
All these approximations may be derived from the explicit Laplace transform provided by the \PK formula \eqref{PK}.

Producing higher order approximations in the  "intermediate regime" when a finite,  but  larger number of moments is known, seems a  very challenging problem.

\beQ In view of the scarcity of
approximations fitting  more than three moments, it is natural to ask {\bf what are the difficulties blocking the development of high  order  moments based approximations}?
\eeQ

Below, we examine
this question in the context of  first passage theory for {\bf spectrally negative Levy processes}.

 \vskip.1truecm{\bf Pad\'e, two-point Pad\'e and other rational   approximations}. It turns out  that most of the approximations currently used are Pad\'e   approximations of Laplace transforms, and that higher order approximations
 are quite easy to
obtain (for example using   the Mathematica command {\bf PadeApproximant}), since the conversion from moments to a {rational Laplace approximation}
 requires  only solving a linear system.
It is also quite easy to produce Pad\'e approximations with specified limiting
   behavior of the Laplace transform at $\I$, so called {\bf two-point Pad\'e approximations} -- see Example \ref{ex:exp}.

   Other
rational approximations of interest
 are those
   minimizing in least squares
 sense the sum of the coefficients -- see Beylkin and Monz\'on \cite{beylkin2005approximation}.
 Another interesting class are
 "{\bf Tijms approximations}", which    try to  incorporate moments fitting with including the exact \CL asymptotics  \eqref{CLa}
as dominant term --see \cite{willmot1998class}.
These may also be obtained by using  a Pad\'e  approximation with a  prescribed pole.

   \vskip.1truecm{\bf The   admissibility of  Pad\'e    approximations in probability}.
While in principle a great tool due to their easiness of implementation, and their convergence for large
$n,$  Pad\'e   approximations (and variations, like two point Pad\'e) applied to Laplace transforms in probability
have the
drawback of the difficulty to  guarantee   {\bf "admissible inverses"}, i.e. {\bf nonnegative densities}
   and {\bf non-increasing survival functions},   when fitting three moments or more is desired
   (note that for fitting two moments $m_1,m_2$ of a nonnegative random  variable, the admissible Gamma approximation
   $$f(x)\sim \T f(x)=\fr{(\mu x)^{\a -1}}{\G(\a)} \mu e^{- \mu x}, \; \a=\frac{m_1^2}{m_2-m_1^2}, \mu=\frac{m_1}{m
   _2-m_1^2}$$   provides an easy solution).


  Even  ensuring the
{\bf nonnegativity  of combinations of exponentials with fixed given rates} is quite a difficult
question (since this
involves an infinite number  of constraints), still open nowadays,  except for two exponentials
(when nonnegativity of $\T f(0)$ and nonnegativity of the coefficient of the asymptotically dominant exponent are clearly necessary and sufficient), and for three exponentials \cite{dehon1982geometric}.

\begin{ex} \la{e:Bot}
Consider  the example due to Harris \cite{harris1992note} (see also
\cite[Ch. 5.4]{fackrell2003characterization}) \bea f(t)= 2 e^{-t} -
6 e^{-2 t} + 6 e^{-3t}=2 e^{-t} - 3( 2 e^{-2 t}) +  2 ( 3
e^{-3t})\\
\Eq \bar F(t)= 2 e^{-t} - 3 e^{-2 t} +  2
e^{-3t}\eea with {\bf canonical coordinates} $(2,-3,2),$
and Laplace transform

{\small \bea f^*(s)=\frac{2 \left(s^2+2 s+3\right)}{(s+1)
   (s+2) (s+3)} =\fr 23 \fr{6}{(s+1)(s+2)(s+3)}-\fr 13 \fr{6}{(s+2)(s+3)}+\fr 23 \fr{3}{s+3}\eea}


Since the "Coxian coordinates" $(2/3,-1/3,2/3),$ produced by the partial fractions decomposition above   
are not nonnegative, this is not a phase-type distribution of order $3$. 

However, by an admissibility criteria for  {\bf combinations of three negative
exponentials} due to
\cite{dehon1982geometric}\fn[4]{with exponents equal to $i=1,2,3$,
this is
$-\w_2 \leq 2 \sqrt{\w_1 \w_3},$ or, after normalization
$$\w_1 + 1 \leq \w_3 \leq \w_1 + 1 + 2 \sqrt{\w_1 }.$$}
(see also
\cite{fackrell2003characterization,bean2008characterization} for a
criterion which allows colliding exponents), we know  this is a proper density. In fact, it is a
phase-type density of order $4$ -- see \eqref{e:4}
 \eeXa

\vskip.1truecm{\bf Admissibility by phase-type
representations}.
 One  "lucky case" in which admissibility is automatic is when one has obtained somehow
      {\bf any phase-type representation} $PH(\bff \a, A)$ 
      with $\bff \a$
      a {\bf probability vector} and $A$ a {\bf subgenerator} matrix
      (satisfying $A_{ij} \geq 0$ for $i \neq j$ and $A \bff 1 \leq \bff 0$).
       In that case, nonnegativity follows from the
probabilistic interpretation of $f(t)$ as the density of the
absorbtion time of the corresponding Markovian semigroup.

      However,  determining when  a phase-type representation exists is
        again notoriously difficult, the so
      called {\bf positive realization}  problem of systems theory
      \fn[4]{Finding the minimum possible order of such a representation is even harder, and known currently only for three moments fitting representations).
Attesting further to the difficulty of providing admissible phase-type approximations are several interesting recent
approaches, like the recursive {\bf minimal order three moments
fitting} formulas of \cite{bobbio2005matching}, 
and the {\bf EM algorithm} approach
 \cite{willmot2011risk}.}.

In example \ref{e:Bot}, it is possible to show that a phase-type of minimal order $4$ is available,
  by  a  recursive approach  of
decomposing the   admissibility domain as union of higher order
admissibility polytopes associated to {phase-type representations}
-- see  \cite{colm1993triangular}.%
One  phase-type representation is:
\beq \la{e:4} A=\left(\begin{array}{llll}
 -1 & 1 & 0 & 0 \\
 0 & -2 & 2 & 0 \\
 0 & 0 & -3 & 3 \\
 0 & 0 & 0 & -4
\end{array}
\right), \q \a=(1/2,0,0,1/2)\eeq

BUTools http://webspn.hit.bme.hu/~telek/tools/butools/butools.html
only obtains a representation of order $5$,
further illustrating the difficulty of this problem.

\beQ Currently,   no algorithmic approach for  testing admissibility of combinations of more than four given negative exponentials (or four terms M\"untz polynomials) is available \fn[4]{a solution might  however be possible by the approach of Faybusovich \cite[Thm 4,5]{faybusovich2002self}, who offers a general representation of the
Koszul-Vinberg  characteristic function  \cite{seeger2012epigraphical} of the positivity cone generated by any
Chebyshev system, as Pfaffian of a matrix of multiple integrals (and taking logarithm
yields a  "self-concordant barrier" function).}.

 \eeQ

The question of providing admissible  approximations
(with
non-specified rates) fitting more than four given values of a Laplace transform (for example moments) is a priori even more difficult. However, two remarkable exceptions in which this the challenging admissibility problem was solved are the {\bf minimal order three moments fitting admissible
    approximations} of Bobbio, Horvath and Telek
\cite{bobbio2005matching}, obtained by a recursive approach on the order, and the non-minimal ones of Johnson and Taaffe \cite{johnson1989matching}, an outcome of the classical moments theory, which  work for any number of moments.

   \vskip.1truecm{\bf Contents and contributions}.
   This problem was motivated by the desire to provide new admissible matrix exponential approximations in ruin theory.

Necessary ruin theory background, including the \PK transform, is reviewed in Section \ref{s:PK}. Section \ref{s:L} presents the "key characters" in our ruin application, the aggregate loss and
 its moments  -- see \eqref{expexp}, \eqref{Lmts}.

First  order Pad\'e approximations of the \PK transform due to Renyi and DeVylder 
are reviewed in Section \ref{s:Ap}. We also provide  here a new approximation, for the perturbed \CL model -- see Theorem \ref{t:per}. 
The remarkable Johnson-Taaffe approximations are reviewed in Section \ref{s:JT}, Theorem \ref{t:JT}. Using their approach, we may "update" a second order approximation due to
 Ramsay to make it work for arbitrary  claims having three moments -- see Theorem \ref{t:Ram}.

In the case of random sums,    it is possible to apply the Johnson-Taaffe approach
both to the individual summands -- we call this a "Ramsay-type approximation", and directly to the sum -- yielding  "Beekman-Bowers-type
approximations".
Theorem \ref{t:new} in Section \ref{s:comp} provides a  comparison between these
two methods, by comparing their explicit {\bf 3 moments
Johnson-Taaffe orders}.

Section \ref{s:prop}  discusses  two-point Pad\'e approximations --see
Theorem \ref{t:Ram2}, whose admissibility is an open problem for the
moment.  Finally, in Section~\ref{sec:num} we provide some numerical
examples.

\sec{Ruin theory background \la{s:PK}}
{\bf The  perturbed Cram\'er
Lundberg risk process} models the reserves of an insurance company
by:

\begin{equation}\label{CLmod} X(t) = u
+ {\c \,} t- S(t) + \s W(t), \q S(t)=\sum_{k=1}^{N_\l(t)} Z_k
\end{equation}
 used in collective risk theory to describe
the surplus $X=\{X(t), t\ge0\}$ of an insurance company. Here, \BEN
\im $u$ is the initial capital, \im $\c \, t$ represents the premium
income up to time $t$, \im $Z_k$ are i.i.d. positive random
variables representing the claims made, with cumulative distribution
function  and density denoted by $F(x)$ and $f(x),$ and (some)
moments denoted by $m_i, i=1,2,...$,  \im $ N=\{N_\l(t), t\ge0\}$ is
an independent Poisson process with intensity $\l$ modeling the
times at which the claims occur, and \im $ W(t)$ is an infinite variation  spectrally negative perturbation, for example a standard Wiener  motion, and $\sigma>0$  is  a scale parameter.\EEN

Since the  jumps of  $X$ are all negative, the moment generating
function $\E[\te{s X(t)}]$ exists for all $s\ge0$ and $t\ge 0$,
and is log-linear in $t.$  The symbol/Laplace exponent/{cumulant generating function}
$\k(s)$ is defined by \beq \label{eq:psi} \k(s)=\log
\left(\E_0 \te{s X(1)}\right) = \boxed{s\le( \c - \l
\bar{F}^*(s)+\T \k(s) \ri)},\eeq where $\bar{F}^*(s)$ denotes the Laplace
transform of the survival function of the claims, and $s \T \k(s) $ is the symbol of the perturbation $\s W(t)$. For example, in the Wiener case $\T \k(s) = \frac{\s^2}{2} s$.

Let $T $ be the first passage time of a stochastic process $X(t)$
below $0$: { \begin{equation*} T : = \inf\{t\ge0: X(t) <
0\}
.\end{equation*} The objects of interest in classical ruin theory are the
 "finite-time" and "eventual" ruin probabilities

\bea &&\rui(t,u)= P_u [T
 \leq t], 
 \q \boxed{\rui(u)= P_u [T < \infty]=P_0[L > u]},  \la{eq:ru}\eea
 where $L=-\underline{X},$ also called maximal aggregate loss, is the negative of all-time infimum of the process \eqref{CLmod},
 started from $0$.

The ultimate ruin probability $\rui(u)$ for \eqref{CLmod}  is not identically $1$ iff
the Levy drift/profit rate
 \beq \boxed{p:= \c - \l m_1 >0:=\l m_1 \th} \eeq
 is positive, in which case  adding the condition $\lim_{x \to
 \I} \rui(x)=0$  determines it uniquely.

\vskip.1truecm{\bf The Pollaczek-Khinchine formula.} Taking
Laplace transform of the Kolmogorov equation for the {\bf
 ultimate ruin probabilities} of the perturbed
Cram\'er-Lundberg model  yields the 
Pollaczek-Khinchine formula:
 \begin{eqnarray}\la{PK}\rui^*(s)=\frac{1}{s}-
\frac{\k'(0)}{\k(s)}=\frac{\r(1-f_e^*(s))+ \T \k(s)/\c}{s(1- \r f_e^*(s)+\T \k(s)/\c)}:=\frac{1}{s}(1-
\f(s))\end{eqnarray} The first expression, in terms of the symbol
$\k(s)$ of the Levy process involved, is valid  for all spectrally negative L\'{e}vy processes,
-- see for example \cite{kyprianou2006introductory}). The second
emphasizes the fact that the result in the Cram\'er-Lundberg case
depends only on the "equilibrium density" of the claims, defined
by \beq \la{eq:sex} f_e(x):=\bar{F}(x)/m_1,\eeq the estimation of which may be a
convenient starting point. The third expression is equivalent to $ \la{eq:rud}\rui(x)=\int_x^\I \T \p(u) du, \q  x>0$ where  $\T \p(u)$ is the {\bf density of the continuous part of the distribution of the aggregate loss}
$L$
 (while $\rui(x) $ is the  survival function of $L$).
The \PK formula for $\f(s)$
 \beq \la{eq:L} \boxed{ \f(s)=\E e^{-s L}=1- s \rui^*(s)=\frac{\k'(0)}{\k(s)/s}=
 \frac{1-\r}{1-\r f_e^*(s)+ \T \k(s)/\c}}. \eeq

In the case $\s=0$,
a beautiful probabilistic interpretation of the (nonperturbed) \PK formula \eqref{eq:L} was
 discovered (independently)
by Benes, Kendall and Dubordieu, by expanding the denominator into a geometric series:
\begin{eqnarray}\la{PKf}&&
\boxed{\f(s)=\fr{1-\r}{1- \r
\WH f_e(s)}=\sum_{k=0}^\I (1-\r) \r^k \WH f_e^{k}(s)}\; \; \; \end{eqnarray} where $f_e(x)=\bar F(x)/m_1$ is the
equilibrium distribution of $Z_i$. This reveals that  $\f(s)$  is the Laplace
transform of
a geometric sum of convolutions of the equilibrium/stationary distribution,
i.e. $\p(x)=(1-\r) \sum_{n=0}^\I \r^k f_e^{(*,k)}(x)$\fn[4]{Another interpretation of
$\f(s)$ is  as
   Laplace transform of the stationary waiting time of the M/G/1 queue --see for example \cite[Thm VIII.5.7]{asmussen2003applied}.}, which may be visualized by examining the "ladders" of the paths (the amounts by which the process $-\c t + S(t)$ jumps to new maxima).

This is the so called \PK "ladder decomposition"    $L=\sum_{n=1}^N L_n$ of the
 "maximal aggregate loss random variable" $L$, with
 $N$  a geometric r.v.
$\Pr[N=k]=(1-\r) \r^{k}, k=0,1,..., \r=\frac{\l m_1}{\c}$,
representing the number of ladders.

\beR {\bf The ladder decomposition, $\s>0$}.
In the Brownian perturbed case, a  beautiful probabilistic interpretation of the  \PK formula
  \eqref{eq:L} was recently discovered by Dufresne-Gerber and rederived  in an elementary way by  Kella, by rewriting \eqref{eq:L} as:
 \beq \la{eq:Kel} \f(s)=\fr{1-\r}{1- \r f_e^*(s)+ \T \k(s)/\c}=\fr{1}{1+ \T \k(s)/\c} \fr{1-\r}{1- \r f_e^*(s) \fr{1}{1+ \T \k(s)/\c}}\eeq
  reflects  the fact that $L$ is an independent sum of a "first creep at the current infimum" (which in the diffusion case is an exponential of rate $\frac{2 \c}{\s^2}$)  and of  an alternating geometric sum
 of "compound Poisson ladders and further creeps" when $\s>0$
 -- see \cite[Fig 2]{dufresne1991risk}, \cite{kella2011class}.

In the case that the Brownian perturbation $\sigma W(t)$ is replaced
by a {\bf general spectrally negative perturbation} $Y$ with non-zero expectation $E[Y(1)]\ge 0$,
a similar ladder decomposition holds true. Let
$$\kappa^*(s) = s \T \k(s)$$
 denote the Laplace exponent of \, $Y$. Note that
by the Wiener-Hopf factorisation--see e.g. \cite{kyprianou2006introductory}, $\T \k(s)$ is Laplace exponent of the possibly killed downward ladder process
(and $s$ is the Laplace exponent of the up-crossing ladder process),  that is, $\T \k$ takes the form
$$\T \k(s) = c + \delta s + \int_{(0,\infty)} (1-\te{-s x}) \nu(dx),$$
where the killing rate $c$ and the drift $\delta$ are non negative constants and the L\'{e}vy measure $\nu$ satisfies
the integrability condition
$$\int_{(0,\infty)}[ 1\wedge x] \nu(dx) < \infty.$$
Then, the formula \eqref{eq:Kel} still holds,
 providing a decomposition of $L$ as an independent sum of the {\bf increment of the ladder process of the perturbation $Y$ at an independent exp$(c)$- random time}
and the geometric sum of further  such increments and ``compound Poisson ladder height increments''.
In particular, if the perturbation is given by a completely asymmetric stable process, i.e.
$\kappa^*(s) = s^\alpha$, $\alpha\in(1,2)$, then $\kappa_-^*(s)=s^{\alpha-1}$
and we identify $\frac{1}{1 + \kappa_-^*(s)/c}$ as the Laplace transform of
the non-negative random variable $\tilde Y^{\alpha-1}_{e(c)}$,
where $Y^{\alpha-1}$ is a stable subordinator with parameter $\alpha-1$.
\eeR

\beR
 Note that the aggregate loss $L$ is the mixture of a discrete mass of
$\f(\I)=1-\rui(0)=1-\r$ at $0$, and of a  continuous random
variable. Letting $\T \f(s)$ denote the Laplace transform of the
density $\T \p(u)$, note the decomposition \beq \la{e:fn}
&&\f(s)=1-\r+ \r \T \f(s) \; \Eq \; \p(x) = (1-\r) \d_0(x) + \r \T
\p(x)\no  \\&& \T \f(s):=\fr{\f(s)-\f(\I)}
{\f(0)-\f(\I)}=\fr{\f(s)-(1-\r)} {\r}=(1-\r)\fr{ f_e^*(s)}{1- \r
 f_e^*(s)}, \eeq
 where we denoted by $\p(u)$ the inverse Laplace transform of $\f(s)$, given by the generalized function
 $\p(u)= \T \p(u) + (1-\r) \d_0(u).$

The behavior at $\I$  of $\f(s)$ distinguishes between
the nonperturbed ($\s=0$) and perturbed case ($\sigma>0$):
\begin{eqnarray}\la{lim} \lim_{s \to \I}
  \f(s)=\lim_{s \to \I}
  1- s \rui^*(s)=1- \rui(0)=\begin{cases}
  1-\frac{\l m_1}{\c}=1-\r=\frac{p}{\c}, \; & \s=0\\
  0, \; & \s>0\end{cases}\end{eqnarray}

  The Laplace transform
   $\f(s)$
is an essential quantity in the 
theory of L\'{e}vy processes
\eeR

\beR The roots of the denominator $\k(s)=0$ in the \PK formula determine the asymptotic behavior of ultimate ruin probabilities. More generally, an important role is played by the roots of the Cram\'er Lundberg equation \beq\la{eq:CL}
\k(s)=q, q >0.\eeq \eeR

\beR In
practice, the true distribution of the claims (and interarrival
times) is of course unknown, and since the \PK formula requires
this knowledge,  it should be viewed more as a theoretical than a practical
tool. It may be argued that  the most reliable information
available in insurance data is  contained in the first few integer
moments, and thus it seems natural to start building approximations by
fitting moments, or, equivalently, by Pad\'e interpolation of the
Laplace transform at the origin.
\eeR

\section{
  The moments of the aggregate loss \la{s:L}}

     From now on, we will assume a classic Brownian perturbation. As an alternative  to Ramsay's approximation of the ladder distribution, we  may approximate   directly the aggregate loss distribution. One advantage is that the second moment of the aggregate loss satisfies automatically the second order representability constraint $c_v \geq 1/2$. This allows focusing on the third moment 
     constraint.

  Consider   the
expansion
 \begin{eqnarray}\la{mtsexp} &&\rui^*(s) =\frac{
 \n_{2,\s} /2-  \n_3 s/6+...}{p+  \n_{2,\s} s/2- \n_3
s^2/6+...} \Eq \\&&\f(s) =\frac{p}{p+  \n_{2,\s} s/2- \n_3
s^2/6+...}\no \\&&=\frac{\th}{1 + \th- \le(1-(\frac{m_2}{m_1}+ \frac{\s^2}{2 \l m_1})s/2+ \frac{m_3}{m_1} s^2/6-\frac{m_4}{m_1}
s^3/4!+...\ri)}\no\\&&=\frac{1-\r}{1-\r\le(1-(\frac{m_2}{m_1}+
\frac{\s^2}{2 \l m_1})s/2+ \frac{m_3}{m_1} s^2/6-\frac{m_4}{m_1}
s^3/4!+...\ri)}\no
\end{eqnarray}where $\n_i=\l m_i, i=0,1,2,...$ are
the moments of the L\'{e}vy measure, and $ \n_{2,\s}= \n_2+
{\s^2}.$


\beR When $\s=0,$ the  expression in the last parenthesis of
\eqref{mtsexp} $$\WH f(s):=1-\frac{m_2}{2 m_1}s+ \frac{m_3}{3 m_1}
s^2/2-\frac{m_4}{4 m_1} s^3/3!+...=\E [e^{-s L_i} |\{L_i
>0\}]$$  has moments $\T
m_i=\frac{m_{i+1}}{(i+1) m_1},$  identifying $L_i$ as  the famous
equilibrium/stationary excess/ladder variable generated by $Z_i$,  with
density $f_e(x)=\bar{F}(x)/m_1$, and {\bf stationary excess moments} $\T
m_i=\frac{m_{i+1}}{(i+1) m_1}$. \eeR

\beR Let $l_i= \E L^i, \l_i= l_i/i!, i \geq 1$ denote the
  moments and
  "factorially reduced moments" obtained by "normalizing" with respect
  to the
exponential moments of the maximal aggregate loss. These may be easily obtained, either by the  recursion equivalent of the \PK formula:
\bea \th \l_n=\T \mu_n +\sum_{k=1}^{n-1} \T \mu_k \l_{n-k}, \q n \geq 1, \T \mu_k:=\fr{\T m_{k}}{k!}=\fr{\mu_{k+1}}{m_1}\eea
or by expanding \eqref{mtsexp} in power series.
The first factorially reduced moments are:
\begin{eqnarray}\la{Lmts} && \l_1=\E L= \frac{\T{m}_1}{\th}+ \fr{\s^2}{ p}= \frac{\l {m}_2+\s^2}{2 p},  \\&&
 \l_2=\E L^2/2!= \frac{\T{m}_2}{2 \th}+  \l_1^2, 
 \;
   \l_3=\E L^3/3!= \frac{\T{m}_3}{3!\th}
  + 2  \frac{\T{m}_1 \T{m}_2/2}{ \th^2}  +   \frac{\T{m}_1^3}{\th^3}
  \no \\&&
   \l_4=\E L^4/4!= \frac{\T{m}_4}{4! \th}
  +  \frac{ 2\T{m}_1 \T{m}_3/3!  +  (\T{m}_2/2)^2
  }{\th^2}  +  \frac{3 \T{m}_1^2 \T{m}_2/2}{\th^3} +
  \frac{\T{m}_1^4  }{\th^4},\no
  \end{eqnarray}
and the mass of the continuous part is $\l_0:=\Pr [L>0]=\r=\frac{1}{1+\th}.$

\eeR

\beR

When $\s=0$, the factorially reduced moments $\l_k$ admit an
interesting decomposition:
\begin{eqnarray*}\la{expexp} &&\f(s)=(1-\r)\frac{1}{1-\r \WH f_e(s)}=\sum_{k=0}^\I (1-\r) \r^k  (f_e^{*}(s))^k=
\\&&1+
 \sum_{k=1}^\I \l_k (-s)^k =
1+\sum_{k=1}^\I (-s)^k \sum'_{i_1+i_2+...=k} \p_{i_1}\p_{i_2}...
\end{eqnarray*}where $\sum'$ denotes sum over all
decompositions of $k$ as a sum, and  \bea \p_i=\frac{\l m_{i+1}}{p
(i+1)!}=\th^{-1}\frac{\T m_{i}}{ i!}=\th^{-1}\WT {\mu}_i>0,\eea
where $\T m_{i} :=\frac{m_{i+1}}{(i+1)m_1}, \mu_i:=\frac{m_i}{i!}$
and $\T \mu_i:=\frac{\T m_i}{i!}$. \eeR


 \beR The moments of the conditioned continuous r.v. $L|L>0$,  necessary for applying  certain results from the literature, may be obtained by dividing by $\r=\fr{1}{1+ \th}$. \eeR

  \beR Note that the corresponding moments $l_i$ may also be viewed as moments of the
  ruin function:
  \bea  \rui_0&:=&\l_1=\int_0^\I \rui(x) dx=\frac{\T{m}_1}{\th}, \; \rui_1:=\l_2=
   \int_0^\I x \rui(x) dx=\frac{\T{m}_2}{2 \th}+   (\frac{\T{m}_1}{\th})^2,\\
  \rui_{k}&:=& \fr{\int_0^\I x^{k}\rui(x) dx}{k!}=  \l_{k+1},...\eea
   \eeR

 \section{Renyi, De Vylder, and a new simple approximation for ruin probabilities 
 \la{s:Ap}}

{\bf Approximations  of ultimate ruin probabilities}.  The problem of approximating ultimate ruin probabilities $\rui(u)$
for the Cram\'er Lundberg model \eqref{CLmod} using data on the
distribution $F(u)$ of the claims is a classic of applied
probability, dating back before 1900.  Its roots may be traced back
to the Danish mathematician TN. Thiele, who
founded  the first insurance company, Hafnia, who
is also  the {\bf inventor of
cumulants}  and of {\bf Thiele continued fractions} useful for
Laplace transform inversion \cite{cohen2007numerical}.

In this section, we make  the 
observation  that the exponential mixture approximations recalled in the introduction are particular  cases of Pad\'e approximations of  Laplace transforms (a theme already present  in Thiele's research preoccupations). We also provide a new simple
approximation in this vein for the perturbed \CL model in Theorem \ref{t:per}.

    \beR Concerning  our application which involves random sums, we make the observation that   it is possible to use  Pad\'e, Johnson-Taaffe and any other rational approximations of Laplace transforms at three levels: \BEN \im for the density of the claims, based on the estimates of $m_i, i \geq 1$
\im for the equilibrium density of the claims, using the
estimates of the equilibrium moments $\T m_i, i \geq 1,$ and  the profit rate $p $.
The second level is intuitively  superior to the first, since the
equilibrium density is monotonically decreasing, even when the claims
density isn't \fn[4]{For example, this  gives rise to a smaller JT index, as illustrated  in Example \ref{ex:uni}}. We will call this Ramsay type approximation.
  \im for
the aggregate loss density transform $\f(s)$ (or,
equivalently, the ccdf  transform $\rui^*(s)$), using directly the moments
$\l_{i+1}, i \geq 0,$ given in \eqref{Lmts}. This amounts
to {\bf working directly with the Pollaczek-Khinchine formula},
instead of approximating its denominator, and so intuitively,
should be better,   at least under certain conditions.\fn[4]{More precisely, the third method is expected to be  better in the case of light tails claims and heavy traffic, while the second is expected to be  better in the case of heavy tails claims and light traffic. The fact that both methods are
 better sometimes is illustrated in Theorem \ref{t:new}, from the point of view of yielding a smaller JT index. Identifying the domains within which   methods two and three are preferable in "boundary cases"  is not an easy task. Note also that a mixture
 of the two has been also proposed \cite{sakurai2004approximating}.}. We will call this Beekman-Bowers type approximation.
 \EEN
  \eeR

 One  Pad\'e
  approximation we consider here is: \begin{eqnarray}\la{Pad}
  (\f)^{(n)}(s)={\mathcal Pade}_{(m,n)}(\f(s))=
  {\mathcal Pade}_{(m,n)}([\frac{p}{\k(s)/s}]_{N})
\end{eqnarray} where ${\mathcal Pade}_{(m,n)}$ denotes the "classic Pad\'e
approximation" based on the Taylor series around $0$, where $m=n$
or $m=n+1$ for the classic/perturbed Cram\'er Lundberg  process,
respectively, and where $[f(s)]_N$ denotes truncation of a power
series to its first $N$ terms, with $N=m+n+1$  (for "theoretical
models" where an expression for $\f(s)$ is available, we may also
take $N=\I$). The first case is applicable to the classic, and the
second to the perturbed Cram\'er Lundberg  model.
 As mentioned, the motivation of \eqref{Pad} is that
   the classic
 De Vylder approximation,  
 is  precisely
 the one point Pad\'e approximation of $\rui^*(s),$ around $0$, of orders $(n-1,n)$,
 with $n=1$.

 A second class of
Pad\'e
  approximations we experiment with  is: \begin{eqnarray}\la{Pad1}
(\f)^{(n)}(s)={\mathcal Pade}_{((m,n,m_1,n_1))}(\f(s))= {\mathcal
Pade}_{((m,n,m_1,n_1))}([\frac{p}{\k(s)/s}]_N)
\end{eqnarray}where ${\mathcal Pade}_{((m,n,m_1,n_1))}$ denotes a
two point Pad\'e approximation based on the power series around
$0$ and $\I$. These are indispensable when dealing with the perturbed model.

\begin{ex} The simplest Ramsay-type approximation is   the {\bf one moment Renyi  exponential
approximation of the equilibrium density}  (which may also be viewed as a Pad\'e $(0,1)$ approximation of the aggregate loss density, which imposes also  the correct limiting behavior $\lim_{s
\to \I}
      s \rui^*(s)=\rui(0)=\r=1/(1+\th)$ of the Laplace transform  at $\I$).
 This  amounts to
    looking for an approximation of the form
    $$\rui^*(s)\approx
    \frac{\r}{s+b_0} \Eq f_e^*(s)\approx
    \frac{b_0/(1-\r)}{s+b_0/(1-\r)},$$
    where $ f_e^*(s)$ denotes  the Laplace transform of the stationary excess density of the claims \eqref{eq:sex} (the two being related by the \PK formula \eqref{PK}).

      Fitting one moment yields $\frac{ 1-\r}{b_0}={\T m_1},$ where $ \T m_1=\frac {m_2}{2  m_1} $ (the first moment
    of the stationary excess density), and   
    \beq \la{e:Re} \rui^*(s)&\approx&
   \frac{\r}{\frac{ 1-\r}{\T m_1}+ s}, \q \Eq \q
     \rui(x)
    \approx \r  e^{-  x (1-\r)/\T m_1},\eeq
    which is also correct in the limit $\r \to 1$ when $\T m_1$ exists see \cite{kalashnikov1997geometric}, \cite[(31)]{grandell2000simple}.

\la{ex:exp}
\end{ex}

\vskip.1truecm{\bf De Vylder's  exponential approximation}.
 The simplest Beekman-Bowers type approximation is De Vylder's \eqref{e:DVL}, one of the most popular
   approximations for ruin probabilities, due to its simplicity and  asymptotic correctness \cite{grandell2000simple}, despite of  its
   its being expressed  in terms  of only the {  first three moments} of the claims.


\vskip.1truecm{\bf Derivation of De Vylder as a two  moments Pad\'e
approximation}. De Vylder's approximation  was obtained originally by
equating the first three moments  of   the original Cram\'er
Lundberg process with those of   a new process
 with
 { exponential claims}, and with different arrival intensity and premium rate.
We check now, using the moments \eqref{Lmts} of the aggregate loss density,
that De Vylder's
formula coincides with the Pad\'e $(0,1)$ approximation around $0$ of
the  Laplace transform  $\rui^*(s)$ of the ruin probabilities.   We start  by expanding in power series the numerator and denominator of the
Pollaczek-Khinchine formula \eqref{PK}:

\bea &&\rui^*(s)=\fr{ \l m_2 /2- \l m_3 s/6+...}{p+ \l m_2 s/2- \l m_3
s^2/6+...}\approx \fr{a}{s + \a} \\&& \Eq a s (p+ \l m_2 s/2- \l m_3
s^2/6+...) \\&&\approx (s+ \a) (\l m_2 s/2- \l m_3 s^2/6+...)\\&&\Eq
 \begin{cases} a p= \a \l m_2 /2,\\  a  m_2 /2= m_2 /2-\a m_3 /6 \end{cases} \Eq
 \begin{cases}
 a= \fr{3 \l m_2^2}{3 \l  m_2^2+2 p  m_3},\\
 \a=\fr{6 p m_2}{ 3 \l
m_2^2+2 p  m_3} \end{cases}\eea
and so
\beq \la{e:DVL}
\rui_{DV}^*(s)=\fr{a}{s + \a} =
\fr{ \fr{3 \l m_2^2}{3 \l  m_2^2+2 p  m_3}}{s+\fr{6 p m_2}{ 3 \l
m_2^2+2 p  m_3}}=\frac{3 \l m_2^2}{s(3 \l  m_2^2+2 p  m_3)+6 p m_2}\eeq

 \beQ

  Identify  domains within which   methods two and three  are "preferable" in some sense, in the simplest case of exponential approximations, i.e. {\bf compare the Renyi and De Vylder approximations}. \eeQ

\begin{ex} \la{ex:per} {\bf The Cram\'er Lundberg
process  with Brownian perturbation $\s B(t)$}. The
Laplace
   exponent  is
\beq \la{eq:ps} \kappa(s)=p s + \frac{\k_2 }{2} s^2 + \sum_{k=3}^\I (-s)^{k} \fr{\l m_k}{k!},\eeq
where $\l m_k$ are the moments of the L\'{e}vy measure, and $\k_2 =
\l m_2 + \s^2$.

Now we have two further unknowns of interest, the
probability of "creeping ruin" and that of "ruin by jump". The respective Laplace transforms   satisfy \beq \la{eq:per}
\rui_d^*(s)=\frac{\s^2/2}{\k'(0)} \f(s)=\frac{\s^2/2}{p} \f(s), \q \rui_j^*(s)=\fr 1s- \frac{\f(s)}{s} -\frac{\s^2/2}{p} \f(s) \eeq

The Pad\'e $K=0, L=1$ approximations
 are  unreasonable,  since they cannot satisfy  the boundary conditions $\rui_d(0)=1, \rui_j(0)=0.$
To satisfy those as well as the equation
\beq \la{eq:mper} \rui_d^*(s)=\frac{\s^2/2}{p} \f(s) \Eq 1-s \rui_j^*(s)= (s+ \fr{2 p}{\s^2}) \rui_d^*(s)\eeq
which follows from \eqref{eq:per},
 we must use at least a  $K=1, L=2$ approximation.

 \beT \la{t:per} Consider the exponential approximation
 \begin{eqnarray} \la{eq:pa}
\rui_d(x)
&=&  \frac{a_d - \m_1}{\m_2-\m_1} e^{-\m_1x} + \frac{ \m_2-a_d }{\m_2-\m_1} e^{-\m_2x} \\
\rui_j(x) &=& \frac{a_j}{ \m_2-\m_1 } [e^{-\m_1 x}   - e^{-\m_2 x}],
\end{eqnarray}
for the "creeping ruin" and  "ruin by jump", which satisfy $\rui_d(0)=1, \rui_j(0)=0.$ Then, by fitting the first two moments of the aggregate loss $L$,
one is led to the following admissible approximation:
\beq \la{eq:a} a_d=\fr{3 m_2}{m_3}, \; a_j=a_d \fr{ \l m_2}{\s^2 }=\fr{3 \l m_2^2}{\s^2 m_3}.
\eeq
and $-\m_1$ and $-\m_2$  the negative roots  of $s^2 + (a_d + a_j + \fr{2 p}{\s^2}) s +  a_d \fr{2 p}{\s^2}=0$ (whose  discriminant is non-negative).

\eeT

 {\bf Proof:}  $\rui^*(s),$   $\rui_d^*(s)=\fr{s+ a_d}{s^2 + b_1 s+ b_0}$ are quotients of  monic polynomials (to satisfy $\lim_{s \to \I} s \rui^*_d(s)=\rui_d(0)=1$), with three free coefficients, but the second condition in \eqref{eq:mper} imposes one more condition $b_0=\fr{2 p}{\s^2} a_d$ (so that $\rui_j^*(s)= \fr{a_j}{s^2 + b_1 s + a_d \fr{2 p}{\s^2} }$), leaving only two free coefficients.
 Finally,  fitting the first two coefficients around $0$ of
\beq \la{eq:dt} \f(s)=1- s\rui^*(s)=\fr{2 p}{\s^2} \rui_d^*(s)\approx \fr{2 p}{\s^2} \fr{s+a_d}{s^2 +(a_d + a_j+ \fr{2 p}{\s^2}) s + a_d \fr{2 p}{\s^2} } \eeq
 yields:
 \bea \begin{cases}\fr {\s^2}{2 p} (1+\fr{a_j}{ a_d})=\fr{\l m_2 +\s^2}{2 p} \Eq
 \fr{a_j}{ a_d}=\fr{\l m_2}{ {\s^2}}\\
 \fr {\s^2}{2 p} \fr{a_j}{ a_d^2}=\fr{\l m_3 }{6 p}\end{cases} \eea
 with solution \eqref{eq:a}. Then,   assuming
 $\m_1 < \m_2$, Laplace inversion yields \eqref{eq:pa}.

Moreover, both $\rui_j $ and $\rui_d$ are  admissible. Indeed,
this is obvious for $\rui_j $, since its initial value  and its dominant coefficient $\frac{a_j}{ \m_2-\m_1 }$ are non-negative.

The same is true for $\rui_d $; indeed, we may    check  that its dominant coefficient is non-negative, i.e.  that $\m_1 \leq a_d$, by noting that $s^2 + (a_d + a_j + 2p/\sigma^2) s + a_d 2p/\sigma^2$
is negative at $s=-a_d$ (since
$
(a_d)^2 + (a_d + a_j + 2p/\sigma^2) (-a_d) + a_d 2p/\sigma^2 = -a_da_j < 0
$,
with $a_d$ and $a_j$ being positive).
Therefore, $\m_1 < a_d < \m_2$.

\beR It is easy to check  that this approximation is exact for exponential claims. Indeed, in that case  the density transform in the parametrization \eqref{eq:dt} is:
 \bea
 &&\f(s)=\fr{p}{\k(s)} = \fr{1}{1+ \fr{\sigma^2}{ 2 p} s +  \fr{\l}{p} s \fr{m_1}{s + m_1^{-1}}}
=  \fr{{s + m_1^{-1}}}{({s + m_1^{-1}})(1+ \fr{\sigma^2}{ 2 p} s) +  \fr{\l}{p} s {m_1}}\\&&=  \fr{s + m_1^{-1}}{\fr{\sigma^2}{ 2 p} s^2 +  s(\fr{\l}{p}  {m_1}+ \fr{\sigma^2}{ 2 p} m_1^{-1}+ 1) + m_1^{-1}}=  \fr{ 2 p} {\sigma^2} \fr{s + m_1^{-1}}{ s^2 +  s(\fr{2 \l {m_1}}    {\sigma^2} + m_1^{-1}+ \fr{ 2 p} {\sigma^2} ) + m_1^{-1} \fr{ 2 p} {\sigma^2} }
\eea with $a_d=1/m_1, a_j=2 \l m_1/\s^2$, and
it is  easy to check  that this coincides with our approximation \eqref{eq:a},  in the case of exponential claims. \eeR

\beR \la{r:per} Another admissible approximation exact in the  exponential case, but fitting now only one moment is
$a_j=2 \l m_1/\s^2,$
$a_d=1/ \T m_1=\fr{2 m_1}{m_2}$ (with similar admissibility proof). \eeR

\end{ex}

\sec{Johnson-Taaffe approximations \la{s:JT}}
This elegant  moment fitting method approximates  by {\bf common order Erlang mixtures approximations}\fn[4]{One hint which may explain this choice comes from
discretizing the well-behaved Post-Widder Laplace inversion method \cite{jagerman1978inversion}, which leads
to Erlang mixtures of common order. Note also that pure Erlang distributions
are the { unique extremal points} minimizing the coefficient of variability,
within phase-type distributions representable at order $n$, and that this is a  dense subclass of phase-type distributions.}:
\beq \la{JTrep} \boxed{f(x) \sim \sum _{i=1}^K w_i e_n(\r_i,x):=\int_0^\I e_n(\r,x)  \nu_K(d \r)}\eeq
where $e_n(\r,x)$ is the Erlang density. Finding a discrete measure $\nu_K(d \r)$ with positive coefficients $w_i$ requires   solving
a classic De Prony  system 
with respect to   {\bf Erlang reduced moments} $$\mu_1=\fr{m_1}{n}, \mu_2=\fr{m_1}{n(n+1)},..., \mu_k=\fr{m_k}{n_{(k)}}, \; \;  n_{(k)}=n(n+1)...(n+k-1)$$ obtained dividing by the moments $n_{(k)}$ of the scaled Erlang density of shape $n$.

 \begin{ex} When $n=1,$ one finds the {\bf factorially reduced moments}
 $$\mu_1=\fr{m_1}{1}, \mu_2=\fr{m_1}{2!}, \mu_3=\fr{m_1}{3!}, ...$$
 which intervene in constructing GHE approximations (which have positive weights, if the moments $m_k$ are far enough from the boundaries of the Stieltjes space of moments $\fr{m_2}{m_1^1}\geq 1, \fr{m_3 m_1}{m_2^1}\geq 1$). \end{ex}

 It turns out that   Erlang reducing leads to "positive moment sequences" $\mu_k$ satisfying $\fr{\mu_2}{\mu_1^1}\geq 1, \fr{\mu_3 \mu_1}{\mu_2^1}\geq 1, ....$ (i.e. having positive Hankel matrices of moments up to any desired degree), for $n$ big enough.   Then, the classic Stieltjes method of creating discrete moment fitting measures with {\bf positive weights} may be applied. 


 \beT \la{t:JT} {\bf Johnson-Taaffe three moments 
approximation}. Let  $\bff m=(m_1, m_2, m_3)$ denote the first three moments of a nonnegative r.v. Let \bea
n^*&=&\lceil \max\{ (m_0 m_2/m_1^2-1)^{-1},
 (m_1 m_3 /m_2^2-1)^{-1} -1 \}\rceil \\&
 =&\lceil \max\{ \fr{1}{\H m_2-1},
 \fr{2 \H m_2-\H m_3}{\H m_3-\H m_2}\}\rceil\eea
 denote the {\bf explicit "JT index of degree" $3$} giving the smallest "JT sector" -- see  \cite[Fig 3-7]{bobbio2005matching} -- containing our target
 moments.
  Let $\bff \mu =(\mu_i, i=1,...3)$ denote the Erlang reduced moments of order
 $n^*$, and let
 $$ \boxed{b_0= \mu_2 - \mu_1^2, b_1=\mu_3 - \mu_1 \mu_2, b_2= \mu_3 \mu_1 - \mu_2^2}$$
   denote the coefficients of the denominator of the classic second order Pad\'e approximant, and
   let
   $$ \m_{1,2}=\fr{2 b_0}{b_1 \mp \sqrt{b_1^2 - 4 b_0
 b_2}}, \; \w_1= \fr{1- \m_1 \mu_1}{\m_2-\m_1},
 \w_2=1-\w_1
 $$
 denote its roots and (non-negative) partial fractions coefficients.

 Then, the Erlang reduced moments satisfy $\bff \mu \in \{ b_0 > 0, b_2> 0 \}$, and
 $$f(x)= \sum_{i=1}^2 \w_i e_{n^*}(\m_i, x), \;
 $$
 is a nonnegative, decreasing "common order
 Erlang $n^*$"  
  density fitting  our  three moments.

 \eeT

Using this result, we can "fix" Ramsay's approximation by sums of two exponentials \cite{ramsay1992practical} to make it work for any valid three moments, at the price of using higher (not minimal) order approximations -- see Theorem \ref{t:Ram}.

\beR {\bf In principle, $2 k+ 1$ moments fitting Johnson-Taaffe Pad\'e
approximations} could be obtained by determining numerically higher order JT
indices ensuring the positivity of  corresponding Hankel
determinants. The Hankel determinants that must be taken into
consideration are: \beq &&\mM_k(\mP)=\{ 
m_4(m_2 - m_1^2) + m_2 (m_3 m_1 -m_2^2)
+ m_3
(m_2 m_1 -m_3^2)\geq 0,  \la{e:St} \no\\
&&m_5 m_4\left(m_1 m_3-m_2^2\right)+m_4\left(-m_3^3+2 m_2
   m_4 m_3-m_1 m_4^2\right)\geq 0,...,\} \eeq

   The inequalities above  define implicitly the {\bf JT index of degree $5$}.
\eeR

 \beT \la{t:Ram} {\bf Ramsay updated}  Consider the ruin problem for the classic (nonperturbed) \CL model \eqref{CLmod}.
 \BEN
      \im  Imposing the correct limiting behavior $\lim_{s \to \I}
      s \rui^*(s)=\rui(0)=\r=1/(1+\th)$  at $\I$, the 
      Pad\'e $(1,2)$  approximations  for the equilibrium density and ruin  transforms are:
    \bea&& f_e^*(s)
    \approx \fr N D :=\frac{b_0 + a_1  s }{b_0 + b_1 s +b_2  s^2}, \\&& \rui^*(s)
    \approx \r \fr {D -N} {{D}-\r N}: =\r \frac{b_2 s +b_1-a_1
     }{b_2 s^2  +\T b_1 s+\T b_0 },\eea
    \bea && b_0=\T \mu_2-\T \mu_1^2, \q 
    b_1 = \T \mu_{3} -\T \mu_{2} \T \mu_{1}, \q b_2= \T \mu_{1} \T \mu_{3}-\T \mu_{2}^2
    \\&& a_1= b_1 -\T \mu_1 b_0 = \T \mu_3  - 2 \T \mu_1 \T \mu_2 + \T \mu_1^3
   \\&& \T b_1 s = b_1 -\r a _1, \q   \T b_0 =(1-\r)  b_0,
   \eea
    where $\T \mu_{2}$ are {\bf factorially reduced equilibrium
    moments}.

    \im  For the approximate inverse transform $f_e(x)$  to be hyperexponential (in particular, { nonnegative and completely monotone}), it is sufficient that $b_0 >0, b_2>0$.



   \im  Higher order moments fitting is possible by Erlang approximations of order $n^*$ big enough, where for three moments $n^*=\lceil \max\{\fr {m_1^2}{m_2-m_1^2},\fr
 {m_2^2}{m_1 m_3 -m_2^2} -1 \}\rceil $.
 \EEN
   \eeT

\sec{Comparison between the JT indices of  Ramsay and Beekman-Bowers
type approximations \la{s:comp}}
 In this section we present a simple comparison between these two approaches, based on comparing their  Johnson-Taaffe indices.

 The  moments of the aggregate loss
$L$  satisfy
\begin{eqnarray}\la{sp2} &&\WH{l}_2=\fr{{l}_2}{{{l}_1^2}}=2+ \th
\frac{\T{m}_2}{\T{m}_1^2}=2+ \th \WH{\T m}_2 \geq  2 + \th >2,\\
&&\WH{l}_3=\fr{{l}_3}{{{l}_1 {l}_2}}=3+  \frac{ \frac{\T{m}_3}{\th}
  +  \frac{\T{m}_1 \T{m}_2}{ \th^2}
  }{
  \frac{\T{m}_2 \T{m}_1}{\th^2}+   2(\frac{\T{m}_1}{\th})^3}=
  3+  \frac{ \th \WH{\T{m}}_2(1+ {\th} {\WH{\T{m}}_3)}
  }{2+ \th \WH{\T{m}}_2} \no \end{eqnarray}

 The second moment   satisfies the necessary  inequality  ${\WH{l}_2}
  \geq \fr{n+1}{n}$ \cite{bobbio2005matching} already with $n=1.$
 We only need to investigate  the  necessary inequality  for the third moment,
 which is:
 \bea n &\geq& \boxed{J(L):=\frac{2 \WH{l}_2-\WH{l}_3}
 { \WH{l}_3-\WH{l}_2}}=\frac{\left(\WH{\T{m}}_2 \WH{\T{m}}_3-2
   \WH{\T{m}}_2^2\right) \theta ^2-4
   \WH{\T{m}}_2 \theta -2}{\theta ^2
   \left(\WH{\T{m}}_2^2-\WH{\T{m}}_2
   \WH{\T{m}}_3\right)-2} \Eq \\&  \fr{n+2}{n+1}&
\leq \boxed{\n(L):=\fr{  \WH{l}_3}{  \WH{l}_2}
} =\frac{\theta  \WH{\T{m}}_2 \left(\theta
\WH{\T{m}}_3+4\right)+6}{\left(\theta  \WH{\T{m}}_2+2\right){}^2}
\eea

Here, $J(L)=\frac{2 \WH{l}_2-\WH{l}_3}
 { \WH{l}_3-\WH{l}_2}$ is a "partial JT index", based on the third moment
 admissibility condition, and the "normalized moment" $\n(L)=\fr{  \WH{l}_3}{  \WH{l}_2}$ is a monotone transformation (since $J=\fr{2-\n}{\n-1}$ is a decreasing function in the relevant range $\n \in (1,\I)$), which has been already used in the literature.

 A rough indication of the performance of the Beekman-Bowers and Ramsay methods will be obtained now by checking which
 of $\n(L)$,  $\n(L_i)$ is higher.

 \beT \la{t:new} a) The partial  J index of the aggregate loss is strictly  smaller than $J(L_i)=\fr{2\WH{\T m}_2-\WH{\T m}_3} {\WH{\T m}_3-\WH{\T m}_2}$  iff \bea \theta \WH{\T m}_2 <\frac{3/2 \WH{\T m}_2-
   \WH{\T m}_3}{  \WH{\T m}_3- \WH{\T m}_2} =\frac{3/2 -\n(L_i)}{  \n(L_i)-1} \eea

 In particular, if $\n(L_i)=\fr{\WH{\T m}_3} {\WH{\T m}_2} \geq  \fr 32,$
 in which case $\n(L) <  \fr 32$ holds as well,
 a  two terms exponential mixture distribution matching  the
first three moments exists only for the equilibrium ladders $L_i$.

b) More generally, let $n$ denote the unique integer such that  $\fr {n+2}{n+1} \leq \n(L_i) \leq  \fr {n+1}{n}, n \geq 2,$ is satisfied. If furthermore $ x_1 < \th \WH{\T m}_2 <x_2,$ where $x_{1,2}=\frac{2 \mp \sqrt{2} \sqrt{n^2 +n-a(n^2-1)
   }}{(n +1)a-(n+2)}$, then
a  $n$'th order Erlang mixture distribution matching  the
first three moments exists only for the equilibrium ladders $L_i$.

c) For any $\th,$ there exists a mixture of Erlang $n$ distribution   matching the
first three moments of $L$ if
 $$n \geq \frac{\WH{\T m}_3}{\WH{\T m}_3-\WH{\T m}_2}.$$
 \eeT

 {\bf Proof:} a) Consider the unimodal function
 $$\n(a,x)=\frac{\WH{l}_3}{ \WH{l}_2} =\fr{a x^2 + 4 x +6}{(x+2)^2},$$
 where $a=\frac{\WH{\T m}_3}{\WH{\T m}_2}, x= \th \WH{\T m}_2$,
 which takes values  $$\fr{3 a -2}{2a -1} \leq  \n(a,x) \leq \max[\n(a,\I)=a,\n(a,0)=3/2]$$
 (the lower minimal value is achieved for $x^*=\fr{1}{a-1} \Eq \th^*=\frac{1}{ \WH{\T m}_3-\WH{\T m}_2}$).

b) Fitting $L$ is possible when
\bea &&\fr{a x^2 + 4 x +6}{(x+2)^2} > \fr{n+2}{n+1} \Eq ((n+1) a -(n+2))x^2 - 4 x + 2(n-1) >0. \eea

In the prescribed range of $a,$ the discriminant $\D=2(n^2 +n-a(n^2-1))$ is always positive, and the inequality holds  when $x < x_1$ or when $x > x_2,$ where $x_{1,2}=\frac{2 \mp \sqrt{2} \sqrt{n^2 +n-a(n^2-1)
   }}{(n +1)a-(n+2)}$.

 c) Follows  by minimizing $\n(\th)$
  in $\th,$ which yields $\th^*=\frac{1}{ \WH{\T m}_3-\WH{\T m}_2}$ and
 $\n(\th^*)=\frac{3 \WH{\T m}_3-2 \WH{\T m}_2}{2 \WH{\T m}_3-\WH{\T m}_2}, J(\th^*)=\frac{\WH{\T m}_3}{\WH{\T m}_3-\WH{\T m}_2}$.

\begin{ex} \la{ex:uni} Consider   a three moments fitting example:  the $U[0,1]$ rv,
  with moments $m_1 =\fr 12, m_2 =\fr 13, m_3 =\fr 14, m_4 =\fr 15$.  The Johnson-Taaffe index is $7$, and so this procedure yields $14$ phases, while the Bobbio, Horvath and Telek \cite{bobbio2005matching} and He-Zhang methods
   \cite{he2006spectral}  yield only $9$ phases.

For the equilibrium distribution, the moments are $\T m_1  =\fr 13, \T m_2 =\fr 16, \T m_3 =\fr {1}{10},$ the normalized moments are
$$\WH{\T m}_2= \fr{\T m_2}{\T m_1^2}=\fr 32, \WH{\T m}_3= \fr{\T m_3 }{\T m_2 \T m_1}= \fr 95$$ and the  Johnson-Taaffe index is $4$.  The Bobbio, Horvath and Telek  method
    yields the minimal order of $3$ phases.

Finally consider  the JT approach for a "Beekman-Bowers"  approximation.
 Computing the partial JT index, we find 
$$ \frac{4 \left(9 \theta ^2+30
   \theta +10\right)}{9 \theta
   ^2+40}.$$  By Theorem \ref{t:new} a), this is strictly less than the $4$ required for the Ramsay approach  when $\th \leq 1$. Taking into account
   the integer part, we find out that $\th \leq 1.5$ still yields a J index
   less or equal  to $4$.
   \end{ex}

\begin{figure}
\begin{center}
\includegraphics[scale=0.7]{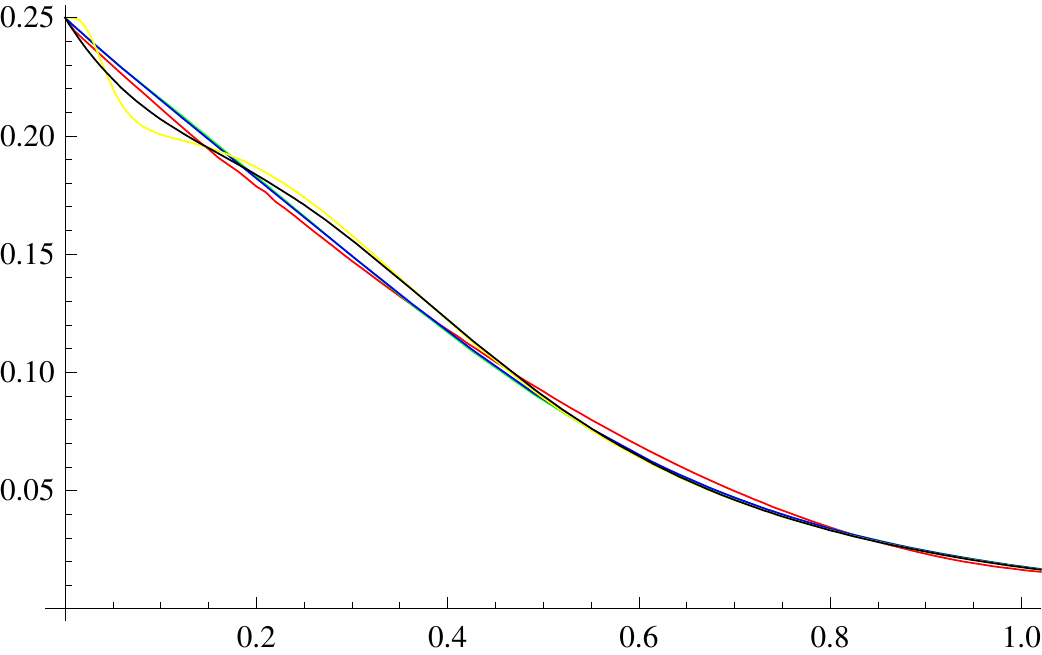}
\end{center}
\caption{Ruin probabilities for  uniform claims.
Red line: original calculated by numerical inverse Laplace;
green line: by 3 moments \cite{johnson1989matching} on the uniform itself;
blue line: by 3 moments \cite{bobbio2005matching} on the uniform itself;
yellow line: by 3 moments \cite{johnson1989matching} on the equilibrium of the uniform;
black line: by 3 moments \cite{bobbio2005matching} on the equilibrium of the uniform.\label{fig:5appr}}
\end{figure}

\sec{Admissible two point Pad{\'e} approximations
\la{s:prop}}

We consider now briefly  including further information beyond moments,  for example 
 by using two-point Pad\'e approximations (of special interest in the case of Levy processes with
 infinite variation paths -- see Example ). Ensuring admissibility in this case is however an open problem.

 We may
formally expand the ruin transform at infinity as well:
\begin{eqnarray}\la{rexp} \rui^*(s)=\int_0^\I e^{-s x} \sum_{k=0}^\I
\rui^{(k)} (0) \frac{ x^k}{k!}=\sum_{k=0}^\I \rui^{(k)} (0)
s^{-k-1}\end{eqnarray} Here,
$$\rui (0)=\begin{cases}1 & \text{if } \s >0\\
\frac{\n_1}{\c}=\r &\text{if } \s =0\end{cases}$$ is well known,
and
 the derivatives at $0$: \begin{eqnarray}\la{der}
 &&\rui' (0)=-\frac{\n_0}{\c}(1-\frac{\n_1}{\c})=- \frac{p \l}{\c^2}=- \frac{\r(1-\r) }{m_1},\\
&&\rui''(0) = -\rui' (0) ( f(0)-\frac{\n_0}{\c} ),
\no\\&&\rui^{(3)}(0) = -\rui' (0) \left(f'(0) +2 \frac{\n_0}{\c}
f(0)  -(\frac{\n_0}{\c} )^2\right),\no\\&& \rui^{(4)}(0) = -\rui'
(0) \left(f''(0)+2 \frac{\n_0}{\c} f'(0)-\frac{\n_0}{\c} f(0)^2 +3
(\frac{\n_0}{\c}) ^2 f(0) -(\frac{\n_0}{\c})
   ^3\right) \no\end{eqnarray} may be obtained recursively, by
differentiating the integro-differential equation 
for
$\rui(x)$.

\beT \la{t:Ram2} {\bf Two-point Pad\'e-Ramsay approximation}
      Imposing both the correct limiting behavior $\lim_{s \to \I}
      s \rui^*(s)=\rui(0)=\r=1/(1+\th)$  at $\I$, and  the first
      derivative at $\I$ $\rui' (0)=- \frac{\r (1-\r)}{m_1}$ leads to the two-point Pad\'e $(1,2)$  approximation:
      \beq \la{eq:twoc} && f_e^*(s)
    \approx \frac{b_0 + a_1  s }{b_0 + b_1 s +b_2  s^2}, \; \; \rui^*(s)
    \approx \r \frac{ b_2 s + b_1   -
    {a_1}  }{\T b_0 +\T b_1 s + b_2 s^2},\eeq
     where $a_1= \fr{b_2}{m_1}  $ and the coefficients $b_i, \T b_i,$  obtained fitting the first two equilibrium/aggregate loss moments, are:
    \bea&& b_2= (2 m_1 m_3-3  m_2^2)/6, \; b_1=
   (m_3-3 m_1 m_2)/3, \; b_0= m_2 -
    2  m_1^2, \\ && \T b_1 =b_1-\r a_1, \q  \T b_0 =(1-\r) b_0.
   \eea

   \eeT

    {\bf Proof:}  The theorem follows from the identity
   \beq \la{eq:id} && s \rui^*(s)= \r \fr{1- f_e^*(s)}{1- \r f_e^*(s)}=\r \frac{b_2  s^2 + (b_1 -a_1)  s }{b_2  s^2 + (b_1 - \r a_1)  s + b_0(1-\r)} \eeq
   by taking into account that we are looking for an approximation of the form
      \bea \rui^*(s) \approx \r \frac{ b_2 s +  b_1   -  \fr{ b_2} {m_1} }{  b_2 s^2+ \T b_1 s+  \T b_0 }:=\r \fr{\T N}{\T D},\eea
      which will satisfy $\lim_{s \to \I} s
   \frac{\T N} {\T D}=1, \lim_{s \to \I} s(s
   \frac{\T N} {\T D}-1) =\rui'(0)/\r= -\frac{1-\r}{m_1}.
   $

\section{Numerical results\label{sec:num}}

\subsection{Without perturbation}

In the following we illustrate the application of the approximations
provided in Theorem~\ref{t:Ram}~and~\ref{t:Ram2}.

\subsubsection{Mixed exponential claim distribution}

With mixed exponential claim distribution
\[
f(x)=
\frac{315 e^{-5 x}}{128} + \frac{7 e^{-4 x}}{8} + \frac{27 e^{-3 x}}{64} +
  \frac{3 e^{-2 x}}{16} + \frac{7 e^{-x}}{128}
\]
and $\lambda=1, c=2/5$ the ruin probability is
\[
\rui(x)=
\frac{245e^{-9x/2}}{32768}
+\frac{135e^{-7x/2}}{8192}
+\frac{567e^{-5x/2}}{16384}
+\frac{735e^{-3x/2}}{8192}
+\frac{19845e^{-x/2}}{32768}
\]

We approximated the ruin probabilities by the  Renyi and De Vylder's first order exponential
approximation given in \eqref{e:Re} and \eqref{e:DVL}, and by the approximations provided in
Theorem~\ref{t:Ram}~and~\ref{t:Ram2}.  In Figure~\ref{fig:mixexpruin} we
show the ruin probabilities with the different approximations and in
Figure~\ref{fig:mixexperror} the relative error with respect to the exact
solution.  As expected, the approximation in Theorem~\ref{t:Ram2}  is better near $x=0,$ and that may be exploited to obtain a better approximation by switching between the formulas in Theorem~\ref{t:Ram2}, Theorem~\ref{t:Ram} when they cross, starting with the first.

\begin{figure}
\begin{center}
\includegraphics[scale=0.6]{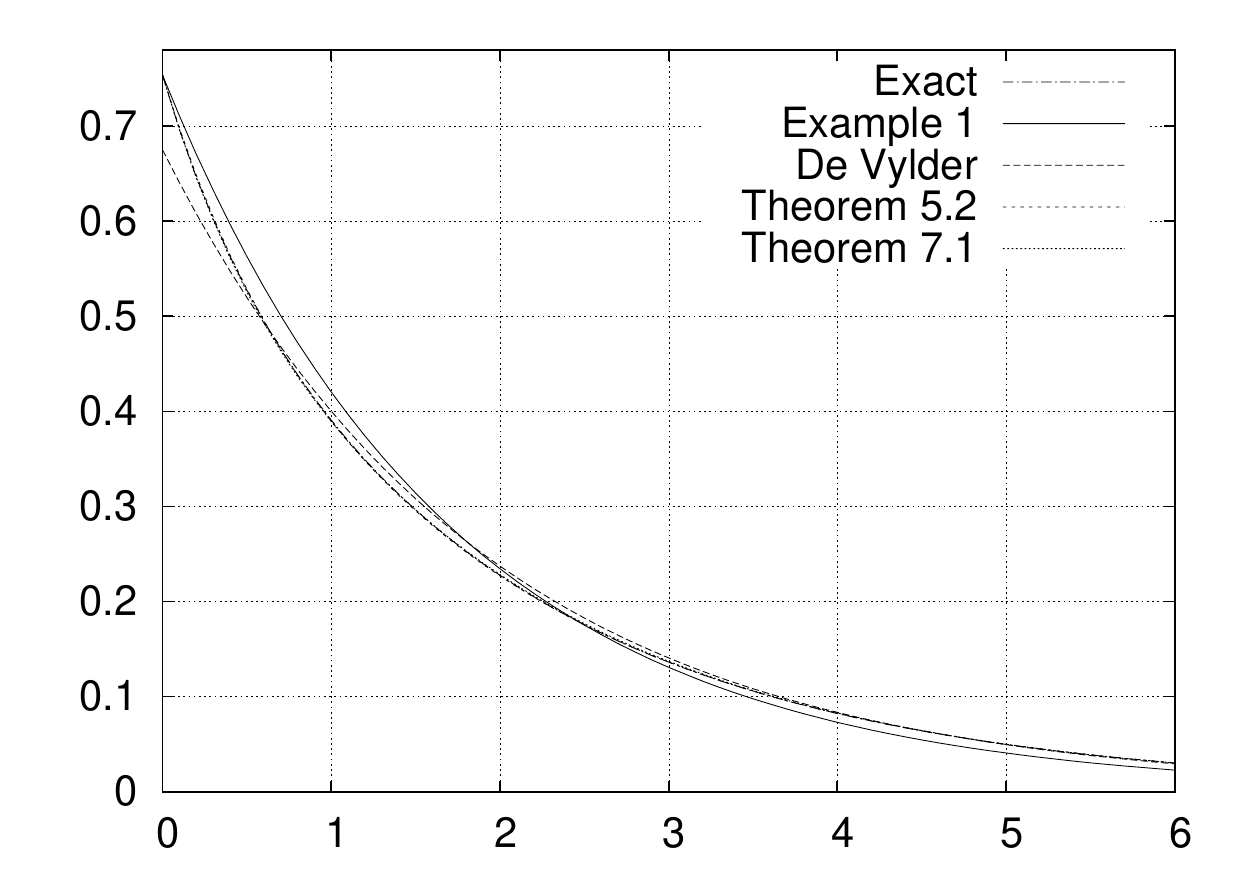}
\end{center}
\caption{Ruin probabilities with mixed exponential claims.\label{fig:mixexpruin}}
\end{figure}

\begin{figure}
\begin{center}
\includegraphics[scale=0.6]{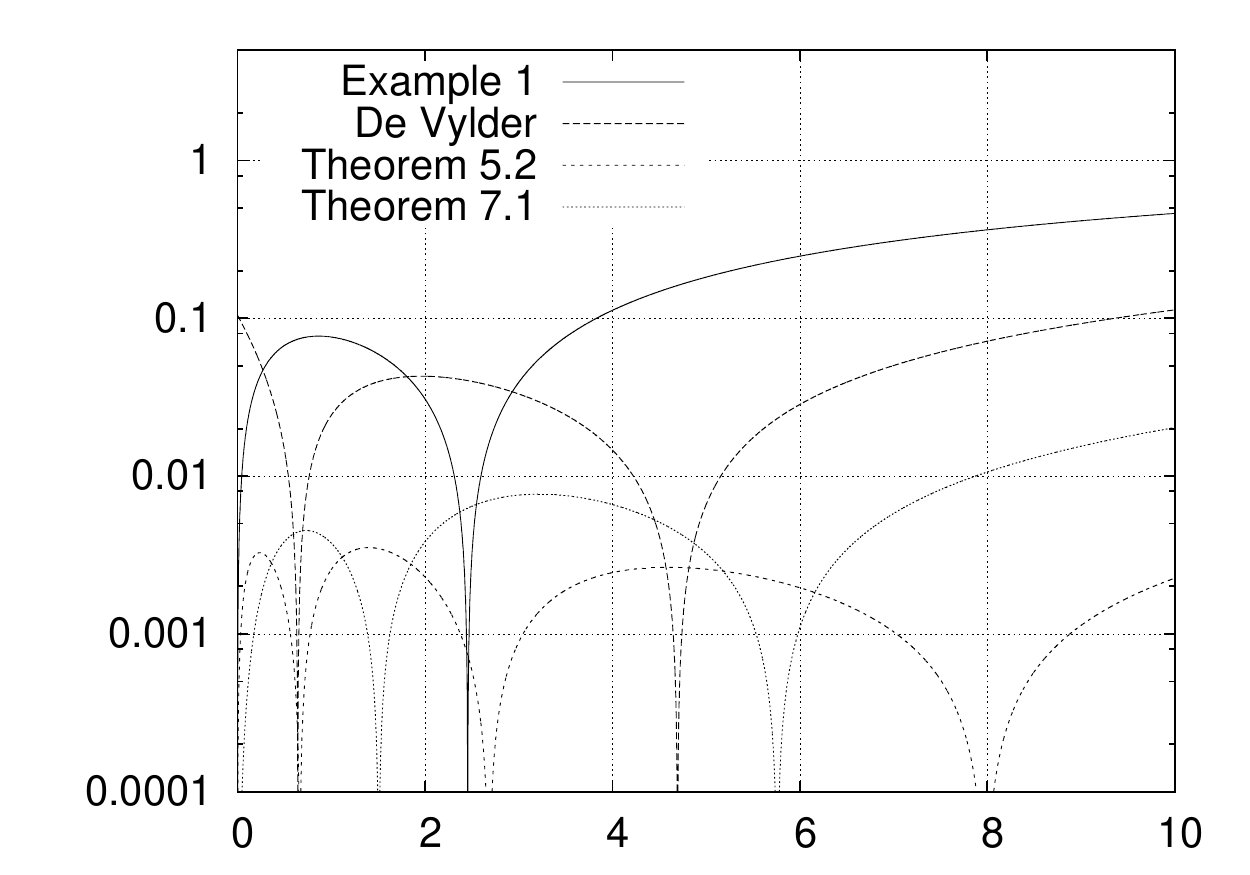}
\end{center}
\caption{Absolute relative error of ruin probabilities with mixed exponential claims.\label{fig:mixexperror}}
\end{figure}

   \subsubsection{Gamma distribution with $\alpha=0.01, \beta=100$}

Let us consider now the Gamma claim distribution
\[
\frac{e^{-\frac{x}{\beta }} x^{\alpha -1} \beta ^{-\alpha }}{\Gamma (\alpha )}
\]
with parameters $\alpha=0.01, \beta=100$ and with claim arrival intensity
$\l=1$ and loading factor $\th=0.1$, which appeared frequently in the
literature (\cite{ramsay1992practical}, \cite{grandell2000simple}).  As for the previous example, we approximated the
ruin probabilities by the formulas given in Renyi,
\eqref{e:DVL} and Theorem~\ref{t:Ram}~and~\ref{t:Ram2}.  Ruin
probabilities themselves are given in Table~\ref{tab:gamma1ruin} while in
Table~\ref{tab:gamma1error} we provide absolute relative errors.

\begin{table}
\begin{center}
\footnotesize
\begin{tabular}{llllll}
\hline
 $x$ & exact~$\rui(x)$ & Renyi & DeVylder &
  Theorem~\ref{t:Ram} & Theorem~\ref{t:Ram2} \\
 \hline
0. & 0.909091 & 0.909091 & 0.882867 & 0.909091 & 0.909091 \\
 300. & 0.521143 & 0.529743 & 0.522539 & 0.521107 & 0.522526 \\
 600. & 0.308668 & 0.30869 & 0.309273 & 0.308713 & 0.309268 \\
 900. & 0.182866 & 0.179879 & 0.183048 & 0.182888 & 0.183047 \\
 1200. & 0.108338 & 0.104818 & 0.10834 & 0.108347 & 0.10834 \\
 1500. & 0.0641841 & 0.0610794 & 0.0641226 & 0.0641869 & 0.0641233 \\
 1800. & 0.0380254 & 0.035592 & 0.037952 & 0.0380257 & 0.0379527 \\
 2100. & 0.0225279 & 0.0207401 & 0.0224625 & 0.0225272 & 0.0224631 \\
 2400. & 0.0133465 & 0.0120856 & 0.0132948 & 0.0133456 & 0.0132953 \\
 2700. & 0.00790706 & 0.00704247 & 0.00786872 & 0.0079062 & 0.00786908 \\
 3000. & 0.00468448 & 0.00410377 & 0.00465722 & 0.0046838 & 0.00465748\\
\hline
\end{tabular}
\caption{Ruin probabilities with Gamma claims ($\alpha=0.01, \beta=100$)
  and its approximations\label{tab:gamma1ruin}}
\end{center}
\end{table}

\begin{table}
\begin{center}
\footnotesize
\begin{tabular}{lllll}
\hline
 $x$ & Renyi & DeVylder &
  Theorem~\ref{t:Ram} & Theorem~\ref{t:Ram2} \\
\hline
 0 & 0 & 0.0288462 & 0 & 0 \\
 300 & 0.0165011 & 0.00267814 & 0.0000688428 & 0.00265325 \\
 600 & 0.0000714085 & 0.0019599 & 0.000146799 & 0.00194387 \\
 900 & 0.0163373 & 0.000993297 & 0.000119325 & 0.000986105 \\
 1200 & 0.0324864 & 0.0000177571 & 0.0000819863 & 0.000019394 \\
 1500 & 0.0483709 & 0.000957418 & 0.0000440626 & 0.00094697 \\
 1800 & 0.0639946 & 0.00193166 & $6.12674\times 10^{-6}$ & 0.00191241 \\
 2100 & 0.0793618 & 0.00290493 & 0.000031798 & 0.00287691 \\
 2400 & 0.0944767 & 0.00387725 & 0.0000697035 & 0.00384047 \\
 2700 & 0.109343 & 0.00484863 & 0.000107635 & 0.00480311 \\
 3000 & 0.123966 & 0.00581907 & 0.000145554 & 0.00576482\\
\hline
\end{tabular}
\caption{Relative error of approximate ruin probabilities with Gamma claims
  ($\alpha=0.01,\beta=100$) \label{tab:gamma1error}}
\end{center}
\end{table}

\subsubsection{Gamma distribution with $\alpha=2.5, \beta=1$}

Next we consider Gamma distributed claims with $\alpha=2.5, \beta=1$ and
$\l=2/5, c=\frac{4}{5} \left(-1+4 \sqrt{2}\right)$.  This case is
interesting as it has been shown in \cite{avram2011moments} that direct,
moment based Pad{\'e} approximation of the claim distribution does not
result in valid distributions. The approximations presented in this paper
leads instead to valid ruin probabilities.  Ruin probabilities themselves
are given in Table~\ref{tab:gamma2ruin} while in
Table~\ref{tab:gamma2error} we provide absolute relative errors.

\begin{table}
\begin{center}
\footnotesize
\begin{tabular}{llllll}
\hline
 $x$ & exact~$\rui(x)$ & Renyi & DeVylder &
  Theorem~\ref{t:Ram} & Theorem~\ref{t:Ram2} \\
 \hline
 0. & 0.268422 & 0.268422 & 0.299749 & 0.268422 & 0.268422 \\
 0.5 & 0.22854 & 0.217791 & 0.237348 & 0.22894 & 0.228126 \\
 1. & 0.189678 & 0.176711 & 0.187938 & 0.189655 & 0.189069 \\
 1.5 & 0.154441 & 0.143379 & 0.148813 & 0.154172 & 0.154016 \\
 2. & 0.124037 & 0.116334 & 0.117834 & 0.123743 & 0.123926 \\
 2.5 & 0.0986589 & 0.0943911 & 0.0933036 & 0.0984496 & 0.0988216 \\
 3. & 0.0779451 & 0.0765868 & 0.07388 & 0.0778418 & 0.0782763 \\
 3.5 & 0.0612929 & 0.0621407 & 0.0584999 & 0.0612758 & 0.0616894 \\
 4. & 0.0480435 & 0.0504196 & 0.0463215 & 0.0480817 & 0.04843 \\
 4.5 & 0.0375759 & 0.0409093 & 0.0366785 & 0.0376414 & 0.0379079 \\
 5. & 0.0293456 & 0.0331929 & 0.0290429 & 0.0294185 & 0.0296037\\
\hline
\end{tabular}
\caption{Ruin probabilities with Gamma claims ($\alpha=2.5, \beta=1$)
  and its approximations\label{tab:gamma2ruin}}
\end{center}
\end{table}

\begin{table}
\begin{center}
\footnotesize
\begin{tabular}{lllll}
\hline
 $x$ & Renyi & DeVylder &
  Theorem~\ref{t:Ram} & Theorem~\ref{t:Ram2} \\
\hline
0. & 0. & 0.116709 & 0. & 0. \\
 0.5 & 0.0470337 & 0.0385392 & 0.00175079 & 0.0018135 \\
 1. & 0.0683674 & 0.00917815 & 0.000125885 & 0.00321361 \\
 1.5 & 0.0716262 & 0.0364389 & 0.00174365 & 0.00275474 \\
 2. & 0.062096 & 0.0500072 & 0.00236637 & 0.000891164 \\
 2.5 & 0.0432581 & 0.0542808 & 0.00212105 & 0.00164952 \\
 3. & 0.0174272 & 0.0521542 & 0.00132515 & 0.00424903 \\
 3.5 & 0.0138333 & 0.045568 & 0.000278168 & 0.00647046 \\
 4. & 0.0494556 & 0.0358424 & 0.000794076 & 0.00804433 \\
 4.5 & 0.0887094 & 0.023884 & 0.00174297 & 0.00883326 \\
 5. & 0.1311 & 0.0103171 & 0.0024838 & 0.00879521 \\
\hline
\end{tabular}
\caption{Relative error of approximate ruin probabilities with Gamma claims
  ($\alpha=2.5,\beta=1$) \label{tab:gamma2error}}
\end{center}
\end{table}

\subsection{With perturbation}

In this section we illustrate the application of the two approximations
given in Theorem~\ref{t:per} and Remark~\ref{r:per}, respectively.

\subsubsection{Mixed exponential claim distribution}

As without perturbation, in this case the exact ruin probabilities can be
calculated by symbolic inversion of the Laplace transform.  We applied the
approximations given in Theorem~\ref{t:per} and Remark~\ref{r:per}  for three values of
$\sigma$, namely, 0.1, 0.5 and 2.  Figures~\ref{fig:rp0.1}-\ref{fig:rp2}
show the exact and approximate ruin probabilities and we depicted the two
components of the ruin probability (by diffusion and by jump) as well.  The two
approximations result in distinguishable curves only in the case
$\sigma=0.1$.  The relative errors for the three values of $\sigma$ are
provided in Figure~\ref{fig:rperr0.1}-\ref{fig:rperr2}.  The error is
smaller for larger values of $\sigma$ and for most values of $x$ the better
approximation of $\rui(x)$ is by the approach of Theorem~\ref{t:per}.

\begin{figure}
\begin{center}
\begin{minipage}{0.495\textwidth}
\includegraphics[width=\textwidth]
{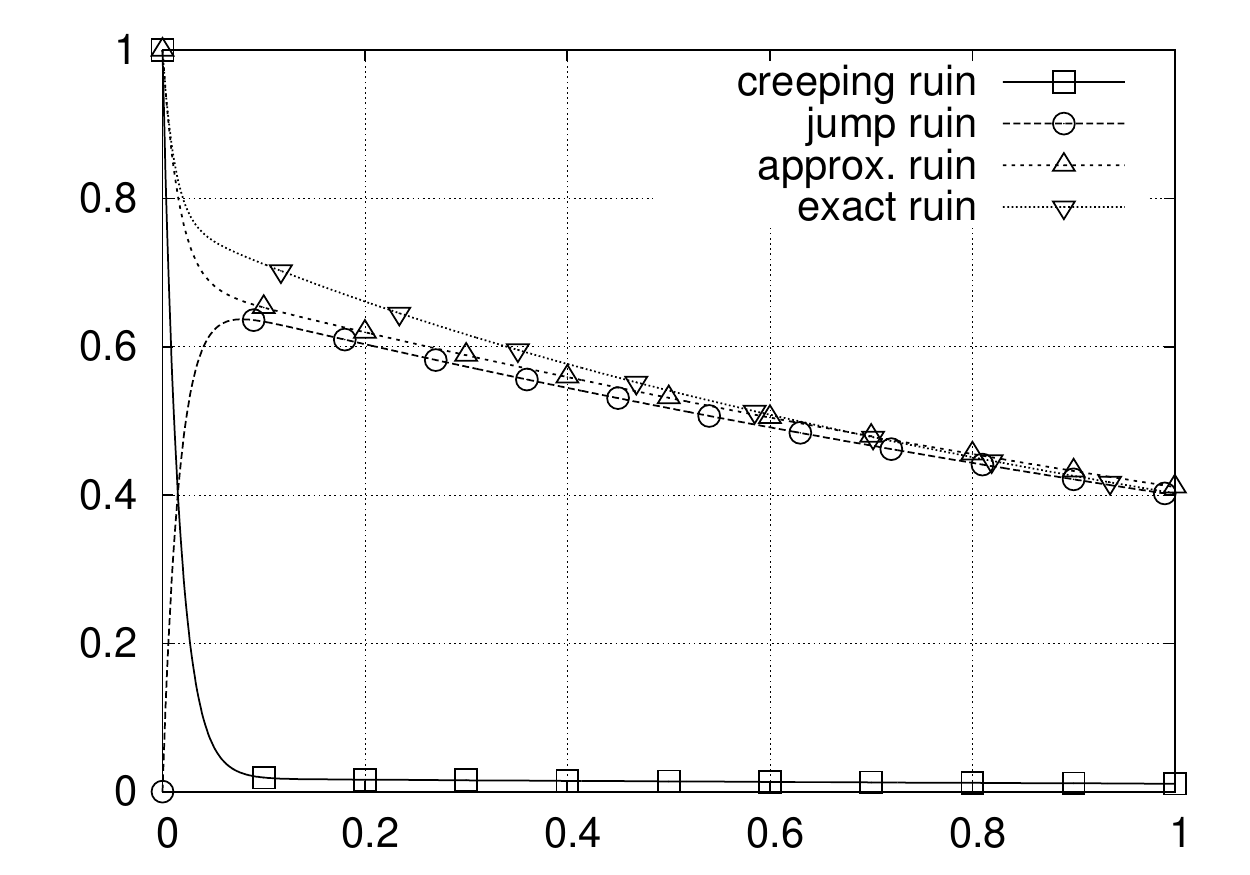}
\end{minipage}
\begin{minipage}{0.495\textwidth}
\includegraphics[width=\textwidth]
{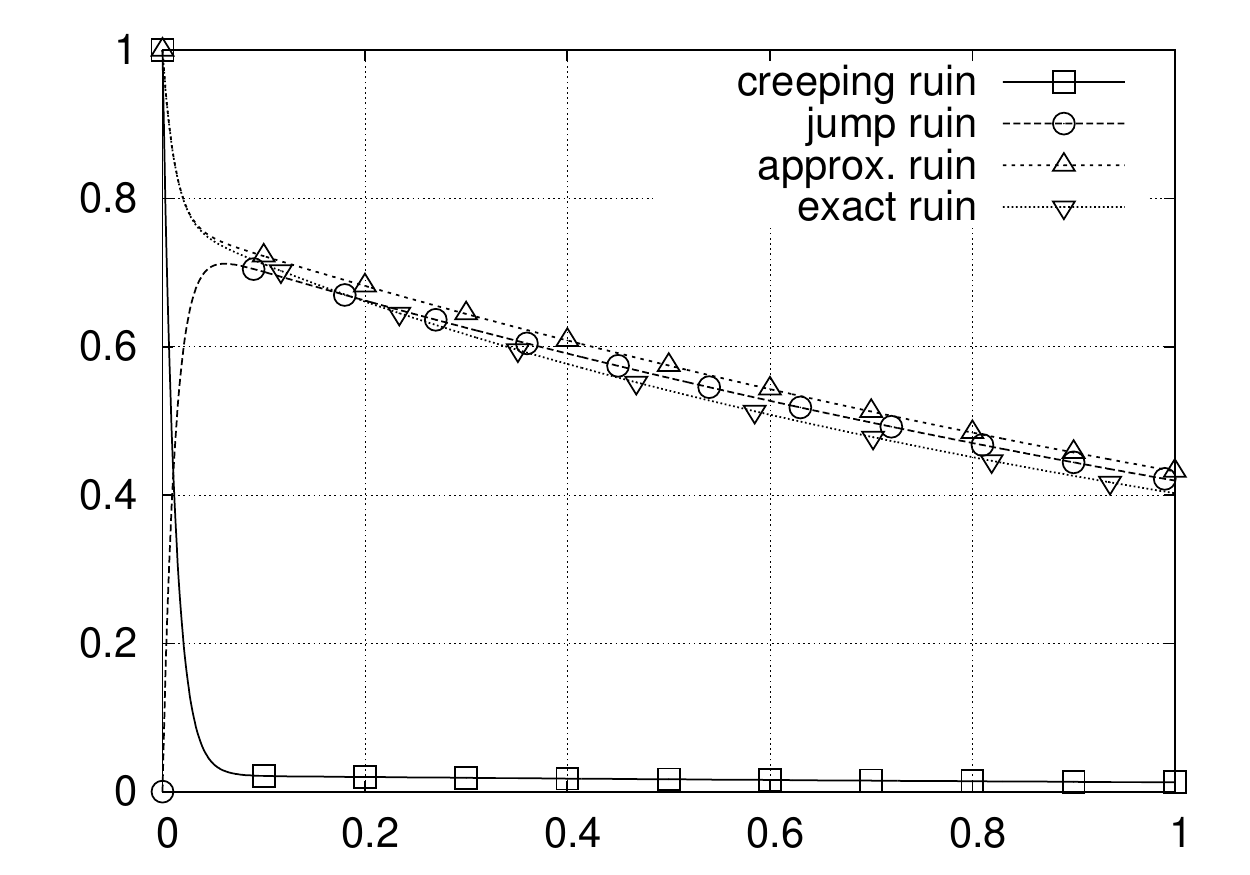}
\end{minipage}
\end{center}
\caption{Ruin probabilities with mixed exponential claim distribution with
  $\sigma=0.1$, with the approximation given in Theorem~\ref{t:per} (left)
  and the one given in Remark~\ref{r:per} (right)
\label{fig:rp0.1}}
\end{figure}

\begin{figure}
\begin{center}
\begin{minipage}{0.495\textwidth}
\includegraphics[width=\textwidth]{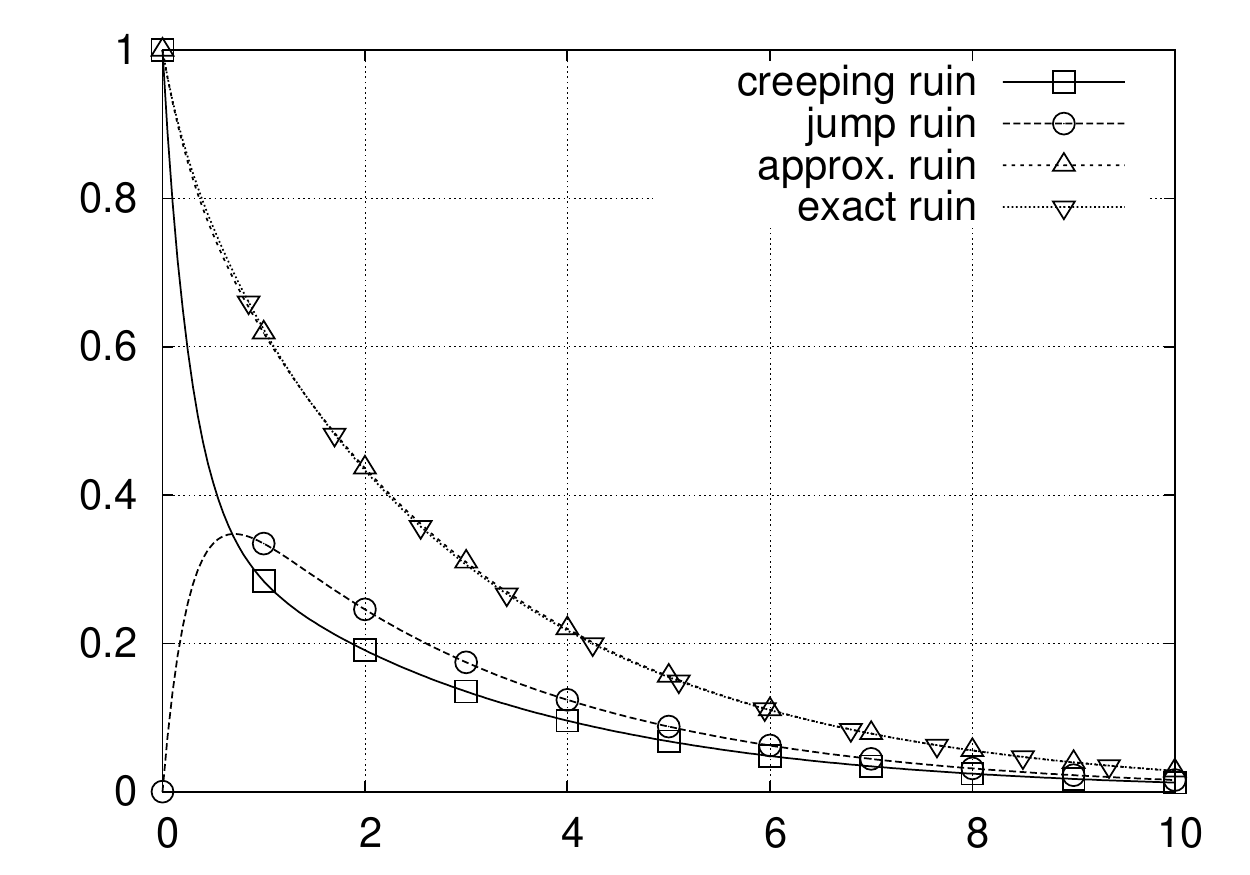}
\end{minipage}
\end{center}
\caption{Ruin probabilities with mixed exponential claim distribution with
  $\sigma=0.5$, with the  approximations
  of Theorem~\ref{t:per} and Remark~\ref{r:per}  (the two approximations look the same)
\label{fig:rp0.5}}
\end{figure}

\begin{figure}
\begin{center}
\begin{minipage}{0.495\textwidth}
\includegraphics[width=\textwidth]{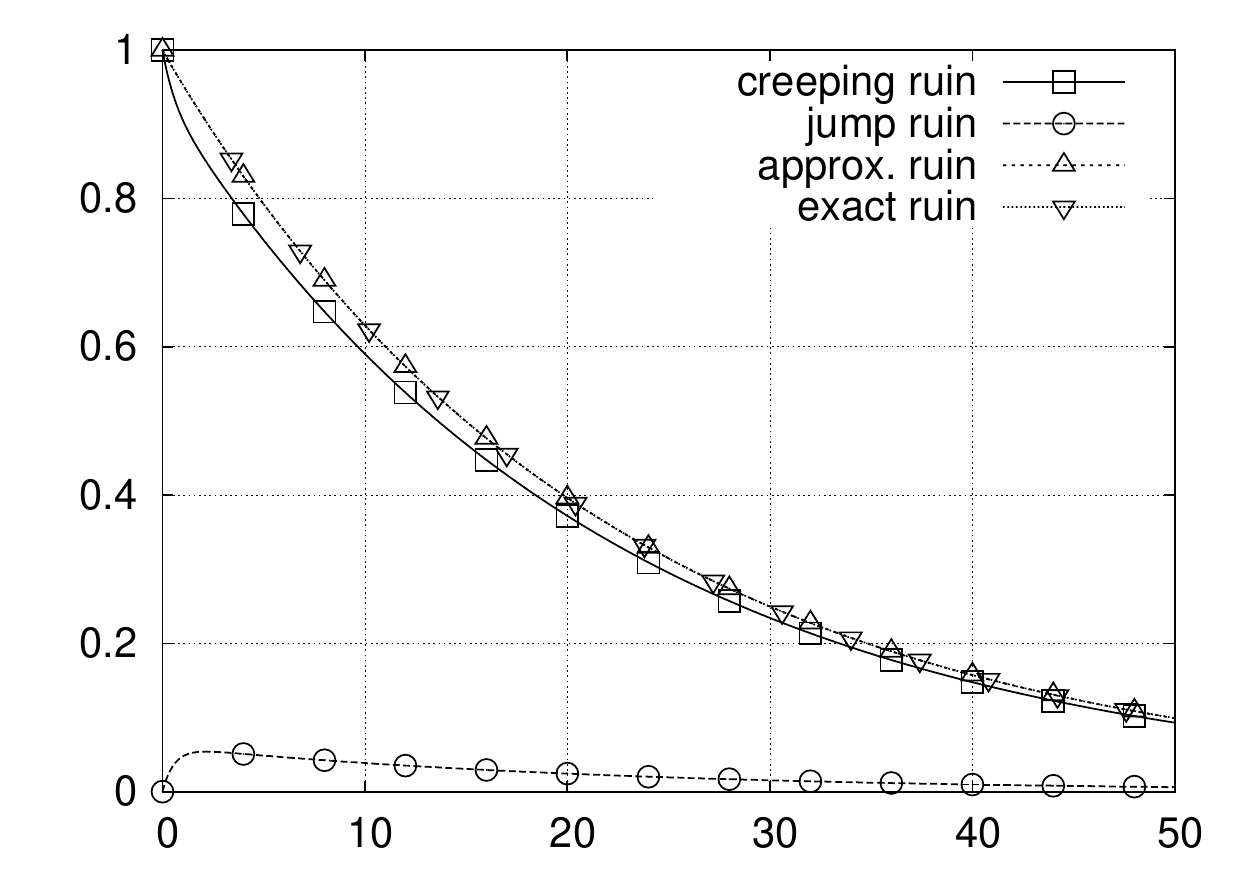}
\end{minipage}
\end{center}
\caption{Ruin probabilities with mixed exponential claim distribution with
  $\sigma=2$, with the  approximations
  of Theorem~\ref{t:per} and Remark~\ref{r:per}  (the two approximations look the same)
\label{fig:rp2}}
\end{figure}

\begin{figure}
\begin{center}
\includegraphics[width=0.495\textwidth]
{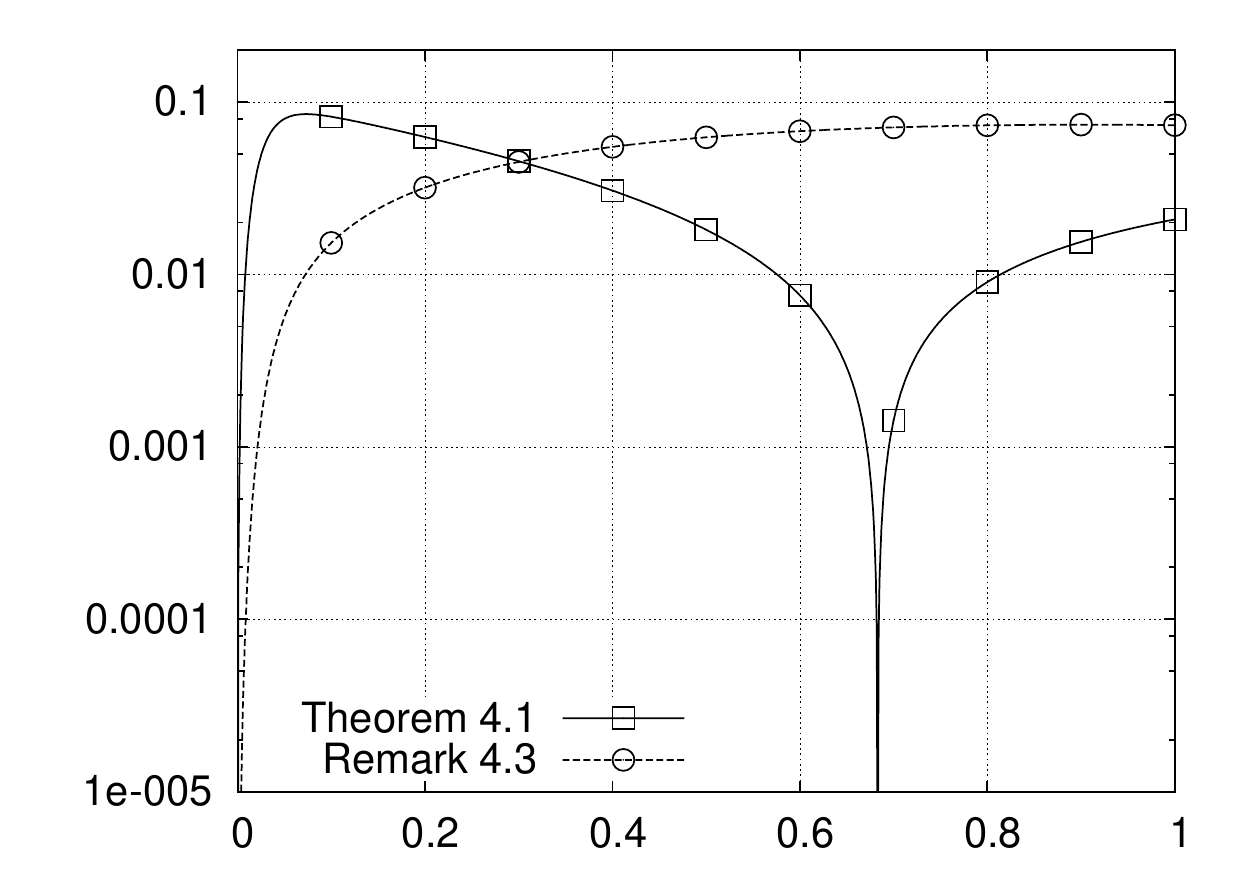}
\end{center}
\caption{Absolute relative error of the approximations given in
  Theorem~\ref{t:per} and Remark~\ref{r:per}  with mixed exponential claims and $\sigma=0.1$.
\label{fig:rperr0.1}}
\end{figure}

\begin{figure}
\begin{center}
\includegraphics[width=0.495\textwidth]
{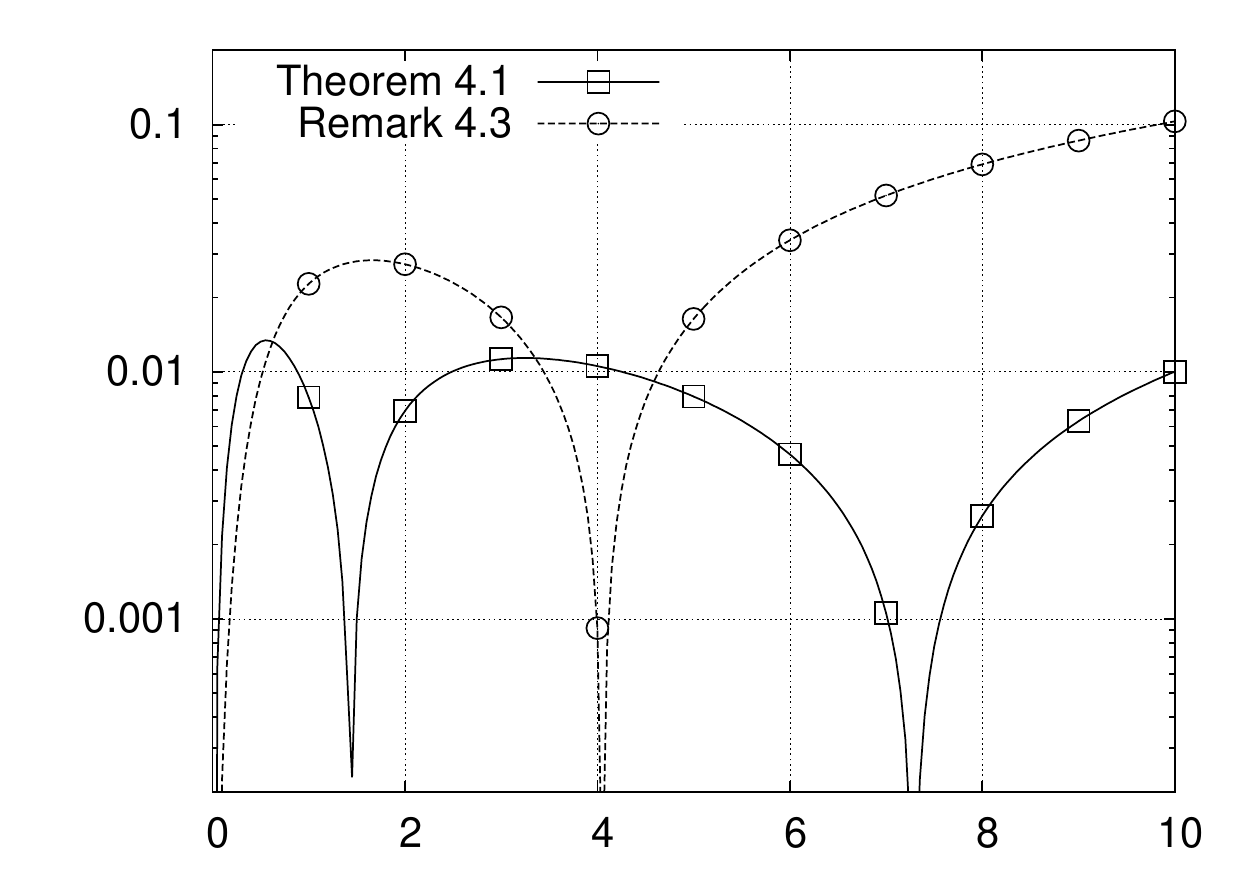}
\end{center}
\caption{Absolute relative error of the approximations given in
  Theorem~\ref{t:per} and Remark~\ref{r:per}  with mixed exponential claims and $\sigma=0.5$.
\label{fig:rperr0.5}}
\end{figure}

\begin{figure}
\begin{center}
\includegraphics[width=0.495\textwidth]
{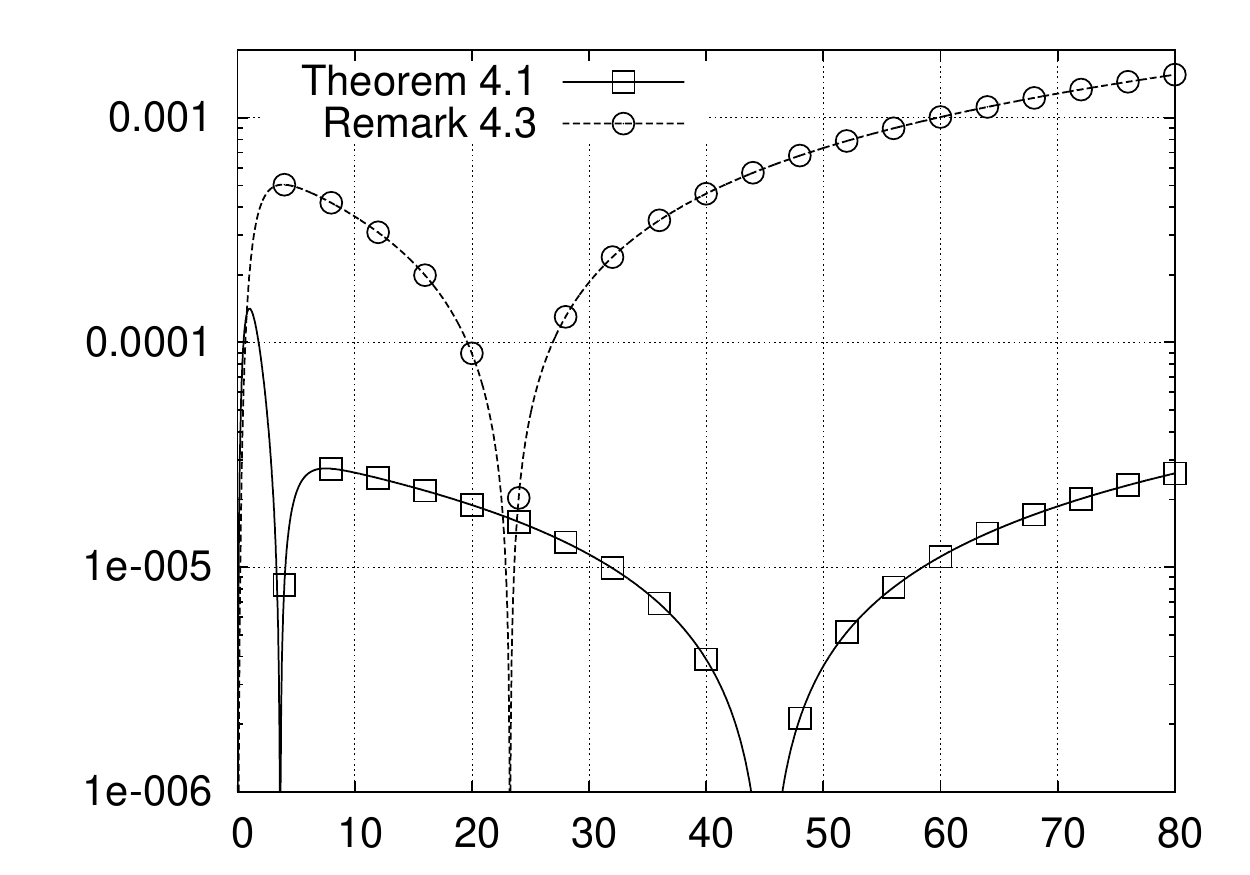}
\end{center}
\caption{Absolute relative error of the approximations given in
  Theorem~\ref{t:per} and Remark~\ref{r:per}  with mixed exponential claims and $\sigma=2$.
\label{fig:rperr2}}
\end{figure}

{\bf Acknowledgement:}  We thank Jiandong Ren for useful comments.

\bibliographystyle{alpha}
\bibliography{Pare}

\end{document}